\documentclass{article}
\usepackage{mathrsfs}
\usepackage{amsmath}
\usepackage{amssymb, amsmath}
\usepackage{graphicx}
\usepackage{fancyhdr}
\usepackage{CJK}
\catcode`\@=11 \@addtoreset{equation}{section}

\catcode`\@=12

\usepackage{color}

\newtheorem{Theorem}{Theorem}[section]
\newtheorem{Proposition}{Proposition}[section]
\newtheorem{Lemma}{Lemma}[section]
\newtheorem{Corollary}{Corollary}[section]
\newtheorem{Definition}{Definition}[section]

\newtheorem{Remark}{Remark}[section]

\newcommand{\newcom}{\newcommand}
\newcommand{\bTheorem}[1]{
\begin{Theorem} \label{T#1} }
\newcommand{\eT}{\end{Theorem}}

\newcommand{\bProposition}[1]{
\begin{Proposition} \label{P#1}}
\newcommand{\eP}{\end{Proposition}}

\newcommand{\bLemma}[1]{
\begin{Lemma} \label{L#1} }
\newcommand{\eL}{\end{Lemma}}

\newcommand{\bCorollary}[1]{
\begin{Corollary} \label{C#1} }
\newcommand{\eC}{\end{Corollary}}

\newcommand{\beq}{\begin{equation}}
\newcommand{\eeq}{\end{equation}}
\newcom{\ben}{\begin{eqnarray}}
\newcom{\een}{\end{eqnarray}}
\newcom{\beno}{\begin{eqnarray*}}
\newcom{\eeno}{\end{eqnarray*}}
\newcom{\bali}{\begin{aligned}}
\newcom{\eali}{\end{aligned}}

\newcommand{\bFormula}[1]{
\begin{equation} \label{#1}}
\newcommand{\eF}{\end{equation}}

\newcommand{\f}{\frac}
\newcommand{\Om}{\Omega}

\newcommand{\p}{\partial}

\newcommand{\vr}{\varrho}

\newcommand{\vt}{\vartheta}
\newcommand{\vu}{\vc{u}}

\newcommand{\vc}[1]{{\boldsymbol #1}}
\newcommand{\Div}{{\rm div}}
\newcommand{\Grad}{\nabla}

\newcommand{\dx}{{\rm d} x}

\newcommand{\dt}{{\rm d} t }

\newcommand{\dxdt}{\dx\dt}

\newcommand{\ep}{\varepsilon}

\font\F=msbm10 scaled 1000
\newcommand{\R}{\mbox{\F R}}

\newcommand\Cbox[2]{%
    \newbox\contentbox%
    \newbox\bkgdbox%
    \setbox\contentbox\hbox to \hsize{%
        \vtop{
            \kern\columnsep
            \hbox to \hsize{%
                \kern\columnsep%
                \advance\hsize by -2\columnsep%
                \setlength{\textwidth}{\hsize}%
                \vbox{
                    \parskip=\baselineskip
                    \parindent=0bp
                    #2
                }%
                \kern\columnsep%
            }%
            \kern\columnsep%
        }%
    }%
    \setbox\bkgdbox\vbox{
        \color{#1}
        \hrule width  \wd\contentbox %
               height \ht\contentbox %
               depth  \dp\contentbox
        \color{black}
    }%
    \wd\bkgdbox=0bp%
    \vbox{\hbox to \hsize{\box\bkgdbox\box\contentbox}}%
    \vskip\baselineskip%
}

\begin{document}


\pagestyle{fancy} \lhead{\color{blue}{Global weak solutions to 2D full non-resistive MHD}} \rhead{\emph{Y.Li and Y.Sun}}

\title{\bf On global-in-time weak solutions to a 2D full compressible non-resistive MHD system}

\author{
Yang Li$^\dag$\,\,\,\,\,\,\,\,\,\,\,\,\,\,\,\,\,\,\,\,\,\,\,Yongzhong Sun$^\ddag$
 \\ \\  $^\dag$School of Mathematical Sciences, \\
Anhui University, Hefei 230601, People's Republic of China \\ Email: lynjum@163.com \\ \\
$^\ddag$Department of Mathematics, \\
Nanjing University, Nanjing 210093, People's Republic of China \\ Email: sunyz@nju.edu.cn \\ \\
}

\maketitle
{ \centerline {\bf Abstract}}

{In this paper, we consider a two-dimensional non-resistive magnetohydrodynamic model, taking the fluctuation of absolute temperature into account. Combining the method of weak convergence developed by Lions \cite{LS}, Feireisl et al. \cite{FN1,FNP} from compressible Navier-Stokes(-Fourier) system and the new technique of variable reduction proposed by Vasseur et al. \cite{VWY} and refined by Novotn\'{y} et al. \cite{NP} from compressible two-fluid models, weak solutions are shown to exist globally in time with finite energy initial data. The result is the  first one on global solvability to full compressible, viscous, non-resistive magnetohydrodynamic system in multi-dimensions with large initial data.
}



\tableofcontents

\section{Introduction}
\subsection{Background and governing equations}\label{back-gov}
Magnetohydrodynamics (``MHD" in short) is concerned with the mutual interactions between electively conducting, heat-conductive fluids and magnetic field. The applications of MHD cover astrophysics, thermonuclear reactions and industry, among many others. Due to its physical importance and mathematical challenges, a lot of efforts have been devoted to theoretical analysis and numerical simulations by mathematicians. In mathematics, a simplified but still physically acceptable full compressible MHD model is described through the following partial differential equations in $(0,T)\times \R^N$ (see \cite{HC,DF}):
\beq\label{fmhd_1}
\p_t \vr +\Div_x (\vr \mathbf{V})=0,
\eeq
\beq\label{fmhd_2}
\p_t (\vr \mathbf{V})+\Div_x (\vr \mathbf{V}\otimes \mathbf{V})+\Grad_x p(\vr,\vt)=\Div_x \mathbb{S}(\vt,\Grad_x \mathbf{V})+\mathbf{curl}_x \mathbf{B}\times \mathbf{B},
\eeq
\beq\label{fmhd_3}
\p_t (\vr s)+\Div_x (\vr s \mathbf{V})+ \Div_x \left( \f{\mathbf{q}}{\vt} \right)
\geq \f{1}{\vt}
\left(
\mathbb{S}(\vt,\Grad_x \mathbf{V}):\Grad_x \mathbf{V}- \f{\mathbf{q} \cdot \Grad_x \vt}{\vt}
+\nu |\mathbf{curl}_x \mathbf{B}|^2
\right),
\eeq
\beq\label{fmhd_4}
\p_t \mathbf{B}=\mathbf{curl}_x (\mathbf{V}\times \mathbf{B})-\mathbf{curl}_x(\nu \mathbf{curl}_x\mathbf{B}),
\eeq
\[
\p_t \left(
\f{1}{2}\vr |\mathbf{V}|^2+\vr e+ \f{1}{2}|\mathbf{B}|^2
\right)+\Div_x \left[
\left(
\f{1}{2}\vr |\mathbf{V}|^2+\vr e+p(\vr,\vt)
\right)\mathbf{V}
\right]
\]
\beq\label{fmhd_5}
=\Div_x \left[
(\mathbf{V}\times \mathbf{B})\times \mathbf{B}+\nu \mathbf{B}\times \mathbf{curl}_x \mathbf{B}-\mathbf{q}+\mathbb{S}(\vt,\Grad_x \mathbf{V}) \mathbf{V}
\right],
\eeq
\beq\label{fmhd_6}
\Div_x \mathbf{B}=0.
\eeq
Here, the unknowns $\vr,\mathbf{V}\in \R^N,\mathbf{B}\in \R^N,\vt$ represent the density, velocity, magnetic field and temperature respectively. $s,e,p$ depending on $\vr,\vt$ are the entropy, internal energy and pressure respectively. $\mathbb{S}(\vt,\Grad_x \mathbf{V})$ is the viscous stress tensor given by
\[
\mathbb{S}(\vt,\Grad_x \mathbf{V})
=\mu (\vt)\left(
\Grad_x \mathbf{V}+\Grad_x^t \mathbf{V}-\f{2}{N}\Div_x  \mathbf{V}\mathbb{I}
\right) +\eta(\vt)\Div_x  \mathbf{V}\mathbb{I},
\]
where $\mu (\vt)$ and $\eta(\vt)$ are the shear and bulk viscosity coefficients respectively. $\mathbf{q}$ is the heat flux obeying Fourier's law
\[
\mathbf{q}=-\kappa(\vt)\Grad_x \vt,
\]
with $\kappa(\vt)$ being the heat conductivity coefficient. $\nu$ is the resistivity coefficient, representing the diffusion of magnetic field. In accordance with \emph{Gibbs' relation},
\beq\label{fmhd_6_1}
\vt Ds(\vr,\vt)=De(\vr,\vt)+p(\vr,\vt)D\left( \f{1}{\vr} \right).
\eeq
Moreover, \emph{thermodynamic stability} conditions read as
\beq\label{fmhd_6_2}
\f{\p p}{\p \vr}(\vr,\vt)>0,\,\, \f{\p e}{\p \vt}(\vr,\vt)>0 \text{    for any  } \vr,\vt>0.
\eeq

Based on the ideas from compressible Navier-Stokes system \cite{FNP,HD1,LS} and Navier-Stokes-Fourier system \cite{FE2}, Ducomet and Feireisl \cite{DF} first proved the existence of global weak solutions to the full compressible MHD system (\ref{fmhd_1})-(\ref{fmhd_6_2}) with finite energy initial data. We refer to Hu and Wang \cite{HW2}, Li and Guo \cite{LG1} for relevant results with general constitutive laws of thermodynamic functions. It should be mentioned that in these works the viscosities are present in the equations of velocity, magnetic field and temperature.

In recent years, a lot of mathematicians focused attention on the incompressible MHD system with partial dissipations. Lin et al. \cite{LXZ} proved global small smooth solutions to 2D viscous non-resistive incompressible MHD system; the case of 3D was later handled by Xu and Zhang \cite{XZ}. We refer to \cite{RWZ,ZT} for more results on 2D viscous non-resistive incompressible MHD system. However, the problem is much more involved when constructing \emph{global-in-time} solutions to the compressible MHD system with partial dissipations. Jiang and Zhang \cite{JZW} first proved global well-posedness for 1D isentropic viscous non-resistive MHD system with large initial data. Global well-posedness to planar full MHD model of viscous non-resistive fluids was shown in \cite{LY2} with large initial data. Again in the simplified 1D regime, existence, uniqueness and stability of weak solutions to viscous non-resistive MHD system were proved in \cite{LY3,LS1D} with large initial data. Relevant results in multi-dimensions seem not so fruitful due to mathematical challenges. Wu and Wu \cite{WW} obtained global small smooth solutions to 2D compressible MHD without magnetic diffusion; Tan and Wang \cite{TW} proved global well-posedness to the corresponding 3D case with small initial data. Li and Sun \cite{LS2D} proved the existence of global weak solutions to a 2D compressible isentropic MHD system of viscous non-resistive fluids with finite energy initial data. Feireisl and Li \cite{FL1} obtained the existence of infinitely many global weak solutions to 3D compressible barotropic MHD system with only magnetic diffusion for any smooth initial data; similar result was also obtained therein for a 2D full MHD system. As a consequence, the existence of \emph{global-in-time} weak solutions to full compressible MHD system of viscous non-resistive fluids in multi-dimensional spaces with \emph{large initial data} remains open.

Restricted by mathematical tools, we consider a simplified situation where the motion of fluids takes places in the plane, while the magnetic field acts on the fluids only in the vertical direction, i.e., by choosing
\begin{equation}\label{fmhd_7}
\left\{\begin{aligned}
& \vr(t,x)=\vr(t,x_1,x_2),\,\, \mathbf{B}=(0,0,b)(t,x_1,x_2), \\
& \vt(t,x)=\vt(t,x_1,x_2),\,\,
\mathbf{V}(t,x)=(u_1,u_2,0)(t,x_1,x_2)=:\vu(t,x_1,x_2).\\
\end{aligned}\right.
\end{equation}
For brevity, we will from now on denote by $x$ the two-dimensional spatial variables $(x_1,x_2)$. Assuming the fluids occupy a bounded smooth domain $\Om \subset \R^2$ and zero resistivity, MHD system (\ref{fmhd_1})-(\ref{fmhd_6}), under (\ref{fmhd_7}), reduces to
\beq\label{fmhd_8}
\p_t \vr +\Div_x (\vr \vu)=0,
\eeq
\beq\label{fmhd_9}
\p_t (\vr \vu)+\Div_x (\vr \vu\otimes \vu)+\Grad_x \left( p(\vr,\vt)+\f{1}{2}b^2\right)=\Div_x \mathbb{S}(\vt,\Grad_x \vu),
\eeq
\beq\label{fmhd_10}
\p_t (\vr s)+\Div_x (\vr s \vu)+ \Div_x \left( \f{\mathbf{q}}{\vt} \right)
\geq \f{1}{\vt}
\left(
\mathbb{S}(\vt,\Grad_x \vu):\Grad_x \vu- \f{\mathbf{q} \cdot \Grad_x \vt}{\vt}
\right),
\eeq
\beq\label{fmhd_11}
\p_t b +\Div_x (b \vu)=0,
\eeq
\[
\p_t \left(
\f{1}{2}\vr |\vu|^2+\vr e+ \f{1}{2}b^2
\right)+\Div_x \left[
\left(
\f{1}{2}\vr |\vu|^2+\vr e+p(\vr,\vt)+b^2
\right)\vu
\right]
\]
\beq\label{fmhd_12}
=\Div_x \left[
 -\mathbf{q}+\mathbb{S}(\vt,\Grad_x \vu) \vu
\right],
\eeq
supplemented with the initial conditions
\begin{equation}\label{fmhd_13}
\left\{\begin{aligned}
& \vr|_{t=0}=\vr_0,\,\, \vr\vu|_{t=0}=\vr_0 \vu_0, \\
& \vr s |_{t=0}=\vr_0 s(\vr_0,\vt_0),\,\,b|_{t=0}=b_0,\\
\end{aligned}\right.
\end{equation}
together with the boundary conditions
\begin{equation}\label{fmhd_14}
\vu|_{\p \Om}=\mathbf{0},\,\, \mathbf{q}\cdot \mathbf{n}|_{\p \Om}=0,
\eeq
where $\mathbf{n}$ is the unit outward normal on $\p \Om$.

\subsection{Structural conditions}\label{struct}

To simplify the presentation and include physically relevant cases, we assume the transport coefficients take the form
\begin{equation}\label{fmhd_15}
\left\{\begin{aligned}
& \mu (\vt)=\mu_0+ \mu_1 \vt,\,\, \mu_0,\mu_1 >0,\,\,\eta=0, \\
& \kappa(\vt)=\kappa_0 +\kappa_2 \vt^2+\kappa_3 \vt^3,\,\, \kappa_0,\kappa_2,\kappa_3>0, \\
\end{aligned}\right.
\end{equation}
in agreement with \cite{FN1}.  To proceed, motivated by the Navier-Stokes-Fourier system \cite{FN1} and the full MHD system \cite{DF}, the pressure $p(\vr,\vt)$ takes the form
\begin{equation}\label{fmhd_16}
p(\vr,\vt)=\underbrace{\vr\vt+\vr^{\gamma}}_{p_M(\vr,\vt)} +\underbrace{\f{a}{3}\vt^4}_{p_R(\vt)}
\eeq
for some $a>0$ and $\gamma>1$. In fact, the term $``p_R(\vt)=\f{a}{3}\vt^4"$ is often referred to the radiation component due to Stefan-Boltzmann law. As a consequence of Gibbs' relation, one finds that the specific internal energy obeys
\begin{equation}\label{fmhd_18}
e(\vr,\vt)=\underbrace{\f{1}{\gamma-1} \vr^{\gamma-1}+c_V \vt }_{e_M(\vr,\vt)}
+\underbrace{a \f{\vt^4}{\vr}}_{ e_R(\vr,\vt)  }.
\end{equation}
Here $c_V$ stands for the specific heat at constant volume, which is assumed to be a positive constant for the sake of simplicity. Finally, we deduce from Gibbs' relation, (\ref{fmhd_16}) and (\ref{fmhd_18}) that
\begin{equation}\label{fmhd_21}
s(\vr,\vt)= \underbrace{\log \f{\vt^{c_V}}{\vr}}_{s_M(\vr,\vt)}+ \underbrace{ \f{4a}{3}\f{\vt^3}{\vr}}
_{s_R(\vr,\vt)}.
\eeq

\begin{Remark}
The term $\vr\vt$ in $p_M(\vr,\vt)$ is usually called Boyle's law. Moreover, $p_M(\vr,\vt)$ with the specific value $\gamma=\f{5}{3}$ may be viewed as a good approximation for liquids, see Chapter 1 in the monograph \cite{FE2} for more discussions with different values of $\gamma$.
\end{Remark}

\section{Weak formulation and main result}\label{ma-re}

\subsection{Definition of weak solutions}
We now introduce the definition of weak solutions.
\begin{Definition}\label{din}
$(\vr,\vu,b,\vt)$ is said to be a weak solution to (\ref{fmhd_8})-(\ref{fmhd_14}) in the time-space domain $(0,T)\times \Om$ provided that
\begin{itemize}
\item { regularities:
\[
\vr\geq 0,\,\, \vt >0 \text{   a.e. in  }(0,T)\times \Om,
\]
\[
\vr \in L^{\infty}(0,T;L^{\gamma}(\Om)) ,\,\,b\in L^{\infty}(0,T;L^2(\Om)),\,\,\vu \in L^2(0,T;W_0^{1,2}(\Om;\R^2)),
\]
\[
\vt\in L^{\infty}(0,T;L^{4}(\Om))\cap L^2(0,T;W^{1,2}(\Om));
\]
}
\item { the continuity equation:
\beq\label{din_1}
\int_0^{T}\int_{\Om} \Big( \vr \p_t \phi +\vr \vu \cdot \Grad_x \phi \Big)\dxdt=
\int_{\Om}
\vr_0 \phi(0,\cdot)  \dx
\eeq
for any $\phi \in C^{\infty}_c([0,T)\times \overline{\Om})$;
    }
\item { balance of momentum:
\[
\int_0^{T}\int_{\Om}
\Big( \vr\vu \cdot \p_t \vc{\phi}+\vr \vu\otimes \vu :\Grad_x \vc{\phi}+
\left(p(\vr,\vt)+\f{1}{2}b^2 \right)\Div_x \vc{\phi}
\]
\beq\label{din_2}
-\mathbb{S}(\vt,\Grad_x \vu):\Grad_x \vc{\phi}\Big ) \dxdt=
\int_{\Om}
\vr_0\vu_0 \cdot \vc{\phi} (0,\cdot)\dx
\eeq
for any $\vc{\phi} \in C_c^{\infty}([0,T)\times \Om;\R^2)$;
  }
\item { principle of entropy production:
\[
\int_0^{T}\int_{\Om}
\left(
\vr s \p_t \phi +\vr s \vu \cdot \Grad_x \phi +\f{\mathbf{q}\cdot \Grad_x \phi}{\vt}
\right) \dxdt
\]
\beq\label{din_3}
+\int_0^{T}\int_{\Om}
\f{\phi}{\vt}\left(
\mathbb{S}(\vt,\Grad_x \vu):\Grad_x \vu- \f{\mathbf{q} \cdot \Grad_x \vt}{\vt}
\right)\dxdt
\leq
\int_{\Om}
\vr_0 s(\vr_0,\vt_0) \phi  (0,\cdot)\dx
\eeq
for any $\phi \in C^{\infty}_c([0,T)\times \overline{\Om})$, $\phi \geq 0$;
  }
\item { non-resistive magnetic equation:
\beq\label{din_4}
\int_0^{T}\int_{\Om} \Big( b \p_t \phi +b \vu \cdot \Grad_x \phi\Big)\dxdt=
\int_{\Om}
b_0 \phi(0,\cdot)\dx
\eeq
for any $\phi \in C_c^{\infty}([0,T)\times \overline{\Om})$;
}
\item { conservation of total energy:
\beq\label{din_5}
\int_{\Om} \left(
\f{1}{2}\vr|\vu|^2+\vr e+ \f{1}{2}b^2
\right)(\tau,x)\dx= \int_{\Om} \left(
\f{1}{2}\vr_0|\vu_0|^2+\vr_0 e_0+ \f{1}{2}b_0^2
\right)\dx
\eeq
for a.e. $\tau \in (0,T)$.
}
\end{itemize}
\end{Definition}
\begin{Remark}
It follows from the regularities of $(\vr,b,\vu)$ and the theory of renormalization due to DiPerna and Lions \cite{DL} that the continuity and Maxwell's equations are satisfied in the renormalized sense:
\[
\p_t f(\vr)+\Div_x (f(\vr)\vu)+(f'(\vr)\vr-f(\vr))\Div_x \vu=0 \,\,\,\, \text{   in  }\mathcal{D}'((0,T)\times \Om),
\]
\[
\p_t f(b)+\Div_x (f(b)\vu)+(f'(b)b-f(b))\Div_x \vu=0\,\,\,\, \text{   in  }\mathcal{D}'((0,T)\times \Om),
\]
for any $f\in C^1(\R)$ such that $f'(z)=0$ if $z$ is sufficiently large.
\end{Remark}

\subsection{Main result}

The main result of this paper is concerned with the existence of global weak solutions with large initial data.
\begin{Theorem}\label{ls_1}
Let the structural hypotheses (\ref{fmhd_15})-(\ref{fmhd_21}) be in force. Let $\gamma>1$ and $\Om\subset\R^2$ be a bounded domain of class $C^{2,\alpha}$, $\alpha>0$. Let the initial data $(\vr_0,\vu_0,b_0,\vt_0)$ be subject to
\begin{equation}\label{din_6}
\left\{\begin{aligned}
& 0<\underline{\vr_0}\leq \vr_0 \leq \overline{\vr_0}<\infty  ,\,\, \text{  a.e. in }\Om, \\
& 0<\underline{b_0}\leq b_0 \leq \overline{b_0}<\infty  ,\,\, \text{  a.e. in }\Om,  \\
&  0<\underline{\vt_0}\leq \vt_0 \leq \overline{\vt_0}<\infty ,\,\, \text{  a.e. in }\Om,  \\
\end{aligned}\right.
\end{equation}
\begin{equation}\label{din_7}
\left\{\begin{aligned}
& \vr_0 s(\vr_0,\vt_0)\in L^1(\Om), \\
&  \vr_0 e(\vr_0,\vt_0)\in L^1(\Om),  \\
&  \vu_0 \in L^2(\Om;\R^2).  \\
\end{aligned}\right.
\end{equation}

Then there exists a global weak solution to the initial-boundary value problem (\ref{fmhd_8})-(\ref{fmhd_14}) in the sense of Definition \ref{din}.
\end{Theorem}

\begin{Remark}
Notice that the conditions (\ref{din_6})$_1$ and (\ref{din_6})$_2$ are designed to ensure the domination between the initial density and magnetic field. Alternatively, one may impose:
\[
\vr_0 \geq 0,\,\,\,\underline{C} \vr_0 \leq b_0 \leq \overline{C} \vr_0\,\, \text{  a.e. in }\Om,
\]
\[
\vr_0 \in L^{\gamma}(\Om),\,\, b_0 \in  L^2(\Om),\, \, \f{ \vc{m}_0}{\sqrt{\vr_0}} \in L^2(\Om;\R^2),
\]
\[
 \vc{m}_0 =\vc{0}\,\text{ a.e. in }    \{x\in \Om: \vr_0(x)=0\},
\]
\[
\vt_0 >0 \text{  a.e. in  }\Om,\,\,\,\,\vr_0 s(\vr_0,\vt_0)\in L^1(\Om),\,\, \vr_0 e(\vr_0,\vt_0)\in L^1(\Om),
\]
for some constants $0< \underline{C}< \overline{C}<\infty$. In a similar manner, the domination condition between two densities has been employed in some previous works of compressible two-fluid models, see \cite{KNC,MMM,NP,VWY}.
\end{Remark}

\begin{Remark}
Our full compressible MHD model (\ref{fmhd_8})-(\ref{fmhd_12}) is reminiscent of the Baer-Nunziato type system for a mixture of two compressible heat-conducting gases considered very recently by Kwon et al. \cite{KNC}. In their paper the weak sequential stability property was proved, leaving the construction of  weak solutions open. Due to the specific structure of our MHD model, a rigorous existence proof is given in this paper.
\end{Remark}

Several key observations on the specific structure of our MHD model and the proof are in order:
\begin{itemize}
\item{
The density $\vr$ and the magnetic field $b$ satisfy the continuity equation simultaneously; whence, at least formally, the quantity $\f{b}{\vr}$ obeys the pure transport equation
\begin{equation}\label{din_10}
\p_t \left(   \f{b}{\vr}  \right) + \vu \cdot \Grad_x \left( \f{b}{\vr}\right)=0.
\eeq
On the one hand, (\ref{din_10}) implies that $\f{b}{\vr}$ is a constant of motion along the particle trajectory trasported by the same velocity field $\vu$. In particular, the domination condition between the initial density and magnetic field is reserved for any $t\in (0,T]$. On the other hand, (\ref{din_10}) enables us to invoke the new technique of variable reduction proposed by Vasseur et al. \cite{VWY} and refined by Novotn\'{y} et al. \cite{NP} (see also \cite{KNC,MMM,WH1}), in the context of compressible two-fluid models. This leads to almost strong convergence of approximate sequences $\left\{\f{b_{\ep}}{\vr_{\ep}}\right\}_{\ep>0}$, see (\ref{exi-117}), which is crucial when passing to the limits.}
\item{
In the process of passing to the limit $n\rightarrow \infty$ in the approximation scheme, see Section \ref{f-g-1}, strong convergences of $\{(\vr_n,b_n)\}_{n\geq 1}$ are indispensable due to its high nonlinearity. Strong convergence of $\{\vr_n\}_{n\geq 1}$ can be carried out in the same manner as Navier-Stokes-Fourier system \cite{FN1}. Even if the density and magnetic field obey the same continuity equation, however, the treatment of $\{b_n\}_{n\geq 1}$ requires more restriction condition. Specifically, we are lacking in the uniform bound of $\{\sqrt{b_n}\vu_n\}_{n\geq 1}$ in $L^{\infty}(0,T;L^2(\Om))$, which is needed when deriving the uniform bound of $ \{\vu_n \cdot \Grad_xb_n\}_{n\geq 1}$ in some $L^p((0,T)\times \Om)$ with $p>1$. This crucial step is finished by combining the domination condition (\ref{exi-27}) and the uniform bound of $\{\sqrt{\vr_n}\vu_n\}_{n\geq 1}$ in $L^{\infty}(0,T;L^2(\Om))$ coming from the basic energy estimate (\ref{exi-19}).
}
\item{
Compared with the Navier-Stokes-Fourier system, the approximate entropy production rate in our MHD model gives rise to new terms, for instance at the level of vanishing artificial viscosity,
\[
\f{\ep}{\vt_{\ep}}|\Grad_x b_{\ep}|^2
+\f{\ep\delta }{2\vt_{\ep}}
\left(\Gamma b_{\ep}^{\Gamma-2}+2\right)|\Grad_x b_{\ep}|^2
.
\]
Fortunately, the two new terms are non-negative, which can be dropped. It plays a crucial role in passing to the limit $\ep\rightarrow 0$, see (\ref{exi-110}).
}
\end{itemize}

The rest of this paper is dedicated to the proof of Theorem \ref{ls_1}. More specifically, In Section \ref{apsch}, we prove global solvability to the approximate regularized problem (\ref{exi-2})-(\ref{exi-6}) by deriving sufficient global a priori estimates, with artificial viscosity $\ep>0$, artificial pressure $\delta>0$. In Section \ref{f-g}, we perform Faedo-Galerkin limit by establishing uniform-in-$n$ estimates and passing to the limits $n\rightarrow \infty$. In Section \ref{vani}, we establish estimates uniformly in $\ep$ and pass to the limits $\ep\rightarrow 0$. Here, the crucial steps are the pointwise convergences of $\{(\vr_{\ep}.\vt_{\ep})\}_{\ep>0}$. Following the same line, we then derive estimates uniformly in $\delta$ and pass to the limits $\delta \rightarrow 0$ in Section \ref{vanpr} with the help of cut-off functions $T_k(\cdot)$ additionally. Finally, possible extensions and some remarks on the main result are discussed in Section \ref{conre}. For the convenience of the reader, several inequalities and lemmas are included as appendix.

\section{Existence of weak solutions}\label{exis}

\subsection{Approximation scheme}\label{apsch}
Assuming that $(\vr,\vu,b,\vt)$ is smooth, we deduce the internal energy equation from that of entropy (\ref{fmhd_10}) and Gibbs' relation (\ref{fmhd_6_1})
\beq\label{exi-1}
\p_t (\vr e)+\Div_x (\vr e \vu) +\Div_x \mathbf{q}=
\mathbb{S}(\vt,\Grad_x \vu):\Grad_x \vu-p(\vr,\vt)\Div_x \vu .
\eeq
Let $X_n\subset C^{2,\alpha}(\overline{\Om};\R^2)$ be a finite-dimensional vector space, of which each element satisfies the non-slip boundary condition and $X_n$ is endowed with the topology induced by $L^2(\Om;\R^2)$-inner product. Moreover, since $\p \Om\in C^{2,\alpha}$, we may further suppose that $\cup_{n=1}^{\infty}X_n$ is dense in $W_0^{1,p}(\Om;\R^2)$ for any $1<p<\infty$. Inspired by the approximate scheme proposed in Navier-Stokes-Fourier system \cite{FE2,FN1}, we consider the following regularized problem, with slight modifications in order to accommodate the presence of magnetic field:
\beq\label{exi-2}
\p_t \vr +\Div_x (\vr \vu)=\ep \Delta \vr,
\eeq
\[
\int_0^{T}\int_{\Om}
\left[ \vr\vu \cdot \p_t \vc{\phi}+\vr \vu\otimes \vu :\Grad_x \vc{\phi}+
\left( p(\vr,\vt)+\f{1}{2}b^2 +\delta \left(\vr^{\Gamma}+\vr^2+b^{\Gamma}+b^2\right)
\right)\Div_x \vc{\phi} \right] \dxdt
\]
\[
=\int_0^{T}\int_{\Om} \Big (
\mathbb{S}(\vt,\Grad_x \vu):\Grad_x \vc{\phi}+ \ep \Grad_x \vr \cdot \Grad_x \vu \cdot \vc{\phi}\Big) \dxdt
\]
\beq\label{exi-3}
-\int_{\Om} \vr_{0,\delta}\vu_{0,\delta}\cdot \vc{\phi} (0,\cdot)    \dx,
\text{     for any    }\vc{\phi} \in C_c^1([0,T);X_n),
\eeq

\[
\p_t (\vr e) + \Div_x (\vr e \vu)-\Div_x\left\{  \left[\kappa(\vt)+\delta \left(\vt^{\Gamma}+\f{1}{\vt}\right)\right]\Grad_x \vt\right\}
\]
\[
=\mathbb{S}(\vt,\Grad_x \vu):\Grad_x \vu -p(\vr,\vt)\Div_x \vu
+\ep |\Grad_x b|^2
\]
\beq\label{exi-4}
+
\ep \delta \left[
\left(\Gamma \vr^{\Gamma-2}+2\right)|\Grad_x \vr|^2
+\left(\Gamma b^{\Gamma-2}+2\right)|\Grad_x b|^2
\right]
+\delta \f{1}{\vt^2}-\ep \vt^{5},
\eeq
\beq\label{exi-5}
\p_t b +\Div_x (b \vu)=\ep \Delta b,
\eeq
supplemented with the boundary conditions ($\mathbf{n}$ is the unit outward normal on $\p \Om$):
\begin{equation}\label{exi-6}
\Grad_x \vr \cdot \mathbf{n}|_{\p \Om}=0,\,\,\vu|_{\p \Om}=\mathbf{0},\,\, \Grad_x \vt \cdot \mathbf{n}|_{\p \Om}=0,\,\,\Grad_x b\cdot \mathbf{n}|_{\p \Om}=0,
\eeq
and the initial conditions:
\begin{equation}\label{exi-7}
(\vr,\vu,b,\vt)|_{t=0}=(\vr_{0,\delta},\vu_{0,\delta},b_{0,\delta},\vt_{0,\delta}).
\eeq
In the approximate problem (\ref{exi-2})-(\ref{exi-6}), $\ep,\delta>0$ are sufficiently small parameters, while $\Gamma>0$ is a sufficiently large parameter; the approximate initial data $\{(\vr_{0,\delta},\vu_{0,\delta},b_{0,\delta},\vt_{0,\delta})\}_{\delta>0}$ are chosen in a suitable way such that
\begin{equation}\label{exi-8}
\left\{\begin{aligned}
& (\vr_{0,\delta},b_{0,\delta},\vt_{0,\delta})\in C^3(\overline{\Om}),\,\,\vu_{0,\delta}
   \in C_c^3(\Om), \\
&  \inf_{\Om}\vr_{0,\delta}>0,\,\, \inf_{\Om} b_{0,\delta}>0,\,\,\inf_{\Om}\vt_{0,\delta}>0,\\
& C_{\ast} \vr_{0,\delta}\leq b_{0,\delta} \leq C^{\ast} \vr_{0,\delta}, \text{   with   }C_{\ast}
   =\inf_{\Om}\f{b_0}{\vr_0},\,\,C^{\ast}
   =\sup_{\Om}\f{b_0}{\vr_0}, \\
& \Grad_x \vr_{0,\delta} \cdot \mathbf{n}|_{\p \Om}=\Grad_x \vt_{0,\delta} \cdot \mathbf{n}|_{\p \Om}=
  \Grad_x b_{0,\delta} \cdot \mathbf{n}|_{\p \Om}=0, \\
&  \vr_{0,\delta} \rightarrow \vr_0 \text{   strongly in   }L^{\gamma}(\Om),\,\,b_{0,\delta} \rightarrow b_0 \text{   strongly in   }L^{2}(\Om), \text{   as   }\delta \rightarrow 0,   \\
&  \sqrt{\vr_{0,\delta}}\vu_{0,\delta}\rightarrow \sqrt{\vr_0}\vu_0 \text{   strongly in   }L^{2}(\Om), \,\,\vt_{0,\delta}\rightarrow \vt_0 \text{   a.e. in   } \Om.   \\
\end{aligned}\right.
\end{equation}

The next lemma shows that the approximate problem (\ref{exi-2})-(\ref{exi-8}) is globally solvable. We only give the sketch of the proof, since it follows largely from compressible Navier-Stokes-Fourier system \cite{FN1} with slight modifications.
\begin{Lemma}\label{sol}
Let $\ep,\delta$ be two positive parameters and the assumptions of Theorem \ref{ls_1} be in force. Then there exists $\Gamma_0>0$ such that for any $\Gamma>\Gamma_0$ the approximate problem (\ref{exi-2})-(\ref{exi-8}) admits a unique strong solution $(\vr,\vu,b,\vt)$ in $(0,T)\times \Om$. Moreover, it holds that
\begin{equation}\label{exi-9}
\left\{\begin{aligned}
& (\vr,b) \in C([0,T];C^{2,\alpha}(\overline{\Om})), \,\, (\p_t \vr,\p_t b)\in C([0,T];C^{0,\alpha}(\overline{\Om})), \\
& \inf_{[0,T]\times \overline{\Om}}\vr >0,\,\, \inf_{[0,T]\times \overline{\Om}} b >0,\,\,
 \vu \in C^1([0,T];X_n),  \\
& C_{\ast}\vr \leq b \leq C^{\ast}\vr \text{    in    }[0,T]\times \overline{\Om}, \\
& \vt \in C([0,T];W^{1,2}(\Om))\cap L^{\infty}((0,T)\times \Om),\,\, \inf_{(0,T)\times \Om}\vt>0,        \\
&  (\p_t \vt,\Delta K_{\delta}(\vt)) \in L^2((0,T)\times \Om) ,\text{    with    }
K_{\delta}(\vt)=\int_1^{\vt} \left[
\kappa(z)+\delta \left(  z^{\Gamma}+\f{1}{z} \right)
\right]  dz   .\\
\end{aligned}\right.
\end{equation}
\end{Lemma}
{\bf{Proof.}} The proof of Lemma \ref{sol} mainly consists of two steps:

\emph{Step 1:} local existence of strong solutions;

\emph{Step 2:} global a priori estimates.

To prove step 1, we first fix the velocity $\vu\in C([0,T];X_n)$ and solve the approximate equations of density, magnetic field and temperature in terms of $\vu$, i.e.,
\[
\vr=\vr[\vu],\,\,b=b[\vu],\,\, \vt=\vt[\vu].
\]
Based on this, we then apply the classical Leray-Schauder fixed point theorem to the approximate momentum equation so as to obtain its local solvability, say on the lifespan $[0,T_{\star}]$. Finally, estimates (\ref{exi-9})$_1$, (\ref{exi-9})$_4$, (\ref{exi-9})$_5$ follow as a direct consequence of maximal regularity theory for linear/quasilinear parabolic equations (see \cite{AM}). We refer to Section 3.4 in \cite{FN1} for similar details, which are nowadays standard arguments.

The rest of the proof is dedicated to deriving global a priori estimates. Obviously, one has
\begin{equation}\label{exi-10}
\int_{\Om}\vr(\tau,x)\dx =\int_{\Om}\vr_{0,\delta}\dx,\,\,
\int_{\Om}b(\tau,x)\dx =\int_{\Om}b_{0,\delta}\dx,\,\,\text{   for any   }\tau \in [0,T_{\star}].
\eeq
Next, we show that the domination condition
\[
C_{\ast} \vr_{0,\delta}\leq b_{0,\delta} \leq C^{\ast} \vr_{0,\delta}, \text{   with   }C_{\ast}
   =\inf_{\Om}\f{b_0}{\vr_0},\,\,C^{\ast}
   =\sup_{\Om}\f{b_0}{\vr_0},
\]
holds for any $t\in [0,T]$. Direct computation shows that
\begin{equation}\label{exi-11}
\left\{\begin{aligned}
& \p_t (b-C^{\ast}\vr)+\Div_x [(b-C^{\ast}\vr)\vu ]=\ep \Delta (b-C^{\ast}\vr)  , \\
& (b-C^{\ast}\vr)|_{t=0}=b_{0,\delta}-C^{\ast}\vr_{0,\delta},  \\
&  \Grad_x (b-C^{\ast}\vr)\cdot \mathbf{n}|_{\p \Om}=0    .  \\
\end{aligned}\right.
\end{equation}
It follows from the maximum principle of parabolic equations that $b-C^{\ast}\vr \leq 0$ in $[0,T]\times \overline{\Om}$. By the same token, $b-C_{\ast}\vr \geq 0$ in $[0,T]\times \overline{\Om}$. Therefore,
\begin{equation}\label{exi-12}
C_{\ast}\vr \leq b \leq C^{\ast}\vr \text{    in    }[0,T]\times \overline{\Om}.
\eeq
To proceed, upon choosing $\chi \vu,\chi \in C_c^1([0,T))$ as test functions in the approximate momentum equation (\ref{exi-3}), we deduce for any $\tau \in [0,T_{\star}]$
\[
\int_{\Om}\left[\f{1}{2}\vr|\vu|^2+\f{1}{2}b^2+
\delta \left( \f{\vr^{\Gamma}  }{\Gamma-1}+\vr^2+ \f{b^{\Gamma}  }{\Gamma-1}+b^2 \right) \right](\tau,x)  \dx
\]
\[
+\int_0^{\tau}\int_{\Om}
\Big(
\mathbb{S}(\vt,\Grad_x \vu):\Grad_x \vu -p(\vr,\vt)\Div_x \vu
\Big)\dxdt
\]
\[
+\int_0^{\tau}\int_{\Om}
\Big\{
\ep |\Grad_x b|^2+
\ep \delta \left[
\left(\Gamma \vr^{\Gamma-2}+2\right)|\Grad_x \vr|^2
+\left(\Gamma b^{\Gamma-2}+2\right)|\Grad_x b|^2
\right]
\Big\} \dxdt
\]
\begin{equation}\label{exi-13}
=
\int_{\Om} \left[\f{1}{2}\vr_{0,\delta}|\vu_{0,\delta}|^2+\f{1}{2}b_{0,\delta}^2+
\delta \left( \f{\vr_{0,\delta}^{\Gamma}  }{\Gamma-1}+\vr_{0,\delta}^2+ \f{b_{0,\delta}^{\Gamma}  }{\Gamma-1}+b_{0,\delta}^2 \right)
\right]  \dx.
\eeq
Integrating the approximate internal energy equation (\ref{exi-4}) over $\Om$ and adding the resulting relation to (\ref{exi-13}) gives rise to the approximate conservation of total energy for any $\tau \in [0,T_{\star}]$
\[
\int_{\Om} \left[\f{1}{2}\vr|\vu|^2+\vr e +\f{1}{2}b^2+
\delta \left( \f{\vr^{\Gamma}  }{\Gamma-1}+\vr^2+ \f{b^{\Gamma}  }{\Gamma-1}+b^2 \right)\right](\tau,x)  \dx
\]
\[
=
\int_{\Om} \Big[   \f{1}{2}\vr_{0,\delta}|\vu_{0,\delta}|^2+\vr_{0,\delta} e(\vr_{0,\delta},\vt_{0,\delta}) +\f{1}{2}b_{0,\delta}^2
\]
\begin{equation}\label{exi-14}
+\delta \left( \f{\vr_{0,\delta}^{\Gamma}  }{\Gamma-1}+\vr_{0,\delta}^2+ \f{b_{0,\delta}^{\Gamma}  }{\Gamma-1}+b_{0,\delta}^2 \right)
\Big] \dx
+\int_0^{\tau}\int_{\Om}
\left(
\delta \f{1}{\vt^2}-\ep \vt^5
\right)\dxdt.
\eeq
Dividing the approximate internal energy equation (\ref{exi-4}) by $\vt$ and making use of Gibbs' relation (\ref{fmhd_6_1}) yields the approximate entropy equation
\[
\p_t[\vr s(\vr,\vt)] +\Div_x [\vr s(\vr,\vt)\vu]-
\Div_x\left\{  \left[\f{\kappa(\vt)}{\vt}+\delta \left(\vt^{\Gamma-1}+\f{1}{\vt^2}\right)\right]\Grad_x \vt\right\}
\]
\[
=\f{1}{\vt}\left\{
\mathbb{S}(\vt,\Grad_x \vu):\Grad_x \vu
+ \left[\f{\kappa(\vt)}{\vt}+\delta \left(\vt^{\Gamma-1}+\f{1}{\vt^2}\right)\right]|\Grad_x \vt|^2
+\delta \f{1}{\vt^2}
\right\}+
\f{\ep}{\vt}|\Grad_x b|^2
\]
\[
+\f{\ep\delta }{\vt} \left[
\left(\Gamma \vr^{\Gamma-2}+2\right)|\Grad_x \vr|^2
+\left(\Gamma b^{\Gamma-2}+2\right)|\Grad_x b|^2
\right]
\]
\begin{equation}\label{exi-15}
+\ep \f{\Delta \vr }{\vt} \left(
\vt s(\vr,\vt) -e(\vr,\vt)-\f{p(\vr,\vt)}{\vr}
\right)
-\ep\vt^4,
\eeq
For any given constant $\overline{\vt}>0$, we first multiply (\ref{exi-15}) by $-\overline{\vt}$, then integrate the resulting relation over $\Om$, finally add the integral to (\ref{exi-14}), to find that
\[
\int_{\Om} \left[\f{1}{2}\vr|\vu|^2+\mathcal{H}_{\overline{\vt}}(\vr,\vt)+\f{1}{2}b^2
+
\delta \left( \f{\vr^{\Gamma}  }{\Gamma-1}+\vr^2+ \f{b^{\Gamma}  }{\Gamma-1}+b^2 \right)\right](\tau,x)  \dx
\]
\[
+\overline{\vt}\int_0^{\tau}\int_{\Om}
\f{1}{\vt}\Big\{
\mathbb{S}(\vt,\Grad_x \vu):\Grad_x \vu
+ \left[\f{\kappa(\vt)}{\vt}+\delta \left(\vt^{\Gamma-1}+\f{1}{\vt^2}\right)\right]|\Grad_x \vt|^2
+\delta \f{1}{\vt^2}
\]
\[
+
\ep|\Grad_x b|^2
+\ep\delta \left[
\left(\Gamma \vr^{\Gamma-2}+2\right)|\Grad_x \vr|^2
+\left(\Gamma b^{\Gamma-2}+2\right)|\Grad_x b|^2
\right]
\Big\} \dxdt
+\int_0^{\tau}\int_{\Om} \ep\vt^5\dxdt
\]
\[
=\int_{\Om} \left[\f{1}{2}\vr_{0,\delta}|\vu_{0,\delta}|^2+\mathcal{H}_{\overline{\vt}}
(\vr_{0,\delta},\vt_{0,\delta})
+\f{1}{2}b_{0,\delta}^2
+
\delta \left( \f{\vr_{0,\delta}^{\Gamma}  }{\Gamma-1}+\vr_{0,\delta}^2+ \f{b_{0,\delta}^{\Gamma}  }{\Gamma-1}+b_{0,\delta}^2 \right)\right] \dx
\]
\begin{equation}\label{exi-16}
+\int_0^{\tau}\int_{\Om} \left( \f{\delta}{\vt^2} +\ep \overline{\vt}\vt^4\right)  \dxdt
-\ep \overline{\vt} \int_0^{\tau}\int_{\Om}
\f{\Delta \vr }{\vt} \left(
\vt s(\vr,\vt) -e(\vr,\vt)-\f{p(\vr,\vt)}{\vr}
\right)
\dxdt,
\eeq
where
\[
\mathcal{H}_{\overline{\vt}}(\vr,\vt):=\vr e(\vr,\vt)-\overline{\vt}\vr s(\vr,\vt).
\]
After a tedious but direct calculation, with the help of Gibbs' relation and the structural hypotheses (\ref{fmhd_16}), (\ref{fmhd_18}), (\ref{fmhd_21}), one finally concludes from (\ref{exi-16}) that
\[
\int_{\Om} \left[\f{1}{2}\vr|\vu|^2+\mathcal{H}_{\overline{\vt}}(\vr,\vt)+\f{1}{2}b^2
+
\delta \left( \f{\vr^{\Gamma}  }{\Gamma-1}+\vr^2+ \f{b^{\Gamma}  }{\Gamma-1}+b^2 \right)\right](\tau,x)  \dx
\]
\begin{equation}\label{exi-17}
+\overline{\vt} \int_0^{\tau}\int_{\Om} \sigma_{\ep,\delta}\dxdt +
\int_0^{\tau}\int_{\Om}\ep \vt^5 \dxdt \leq C, \text{   for any   }\tau \in [0,T_{\star}],
\eeq
where $C>0$ depends only on the initial data (\ref{din_6})-(\ref{din_7}), while independent of $n,T_{\star},\ep,\delta$;
\[
\sigma_{\ep,\delta}:=
\f{1}{\vt}\Big\{
\mathbb{S}(\vt,\Grad_x \vu):\Grad_x \vu
+ \left[\f{\kappa(\vt)}{\vt}+\delta \left(\vt^{\Gamma-1}+\f{1}{\vt^2}\right)\right]|\Grad_x \vt|^2
+\delta \f{1}{\vt^2}
\]
\begin{equation}\label{exi-18}
+
\ep|\Grad_x b|^2
+\ep\delta \left[
\left(\Gamma \vr^{\Gamma-2}+2\right)|\Grad_x \vr|^2
+\left(\Gamma b^{\Gamma-2}+2\right)|\Grad_x b|^2
\right]
\Big\}
+\f{\ep}{\vr\vt} \f{\p p_M}{\p \vr}(\vr,\vt)|\Grad_x \vr|^2.
\eeq

Based on the global estimate (\ref{exi-17}), the local solution may be extended within finite steps to a global one. This completes the proof of Lemma \ref{sol}.  $\Box$

\subsection{Faedo-Galerkin limit}\label{f-g}
The aim of this subsection is to perform the Faedo-Galerkin limit $n\rightarrow \infty$ in the approximate problem (\ref{exi-2})-(\ref{exi-8}), with $\ep,\delta$ being positive parameters. To this end, let $\{(\vr_n,\vu_n,b_n,\vt_n)\}_{n=1}^{\infty}$ be a sequence of solutions ensured by Lemma \ref{sol}. Notice that (\ref{exi-12}) and (\ref{exi-17}) give the following uniform estimates:
\begin{equation}\label{exi-19}
\sup_{\tau\in (0,T)}
\int_{\Om} \left[\f{1}{2}\vr_{n}|\vu_{n}|^2+\mathcal{H}_{\overline{\vt}}(\vr_{n},\vt_{n})+\f{1}{2}b_{n}^2
+
\delta \left( \f{\vr_{n}^{\Gamma}  }{\Gamma-1}+\vr_{n}^2+ \f{b_{n}^{\Gamma}  }{\Gamma-1}+b_{n}^2 \right)\right](\tau,x)  \dx\leq C,
\eeq
\begin{equation}\label{exi-20}
\int_0^{T}\int_{\Om}
\f{1}{\vt_{n}}\left\{
\mathbb{S}(\vt_{n},\Grad_x \vu_{n}):\Grad_x \vu_{n}
+ \left[\f{\kappa(\vt_{n})}{\vt_{n}}+\delta \left(\vt_{n}^{\Gamma-1}+\f{1}{\vt_{n}^2}\right)\right]|\Grad_x \vt_{n}|^2
\right\} \dxdt \leq C,
\eeq
\begin{equation}\label{exi-21}
\int_0^{T}\int_{\Om}\left(
\delta \f{1}{\vt_{n}^3}+ \ep\vt_{n}^5+\f{\ep}{\vt_{n}}|\Grad_x b_{n}|^2\right)
\dxdt\leq C,
\eeq
\begin{equation}\label{exi-22}
\int_0^{T}\int_{\Om}\left(
\f{\ep\delta}{\vt_{n}} \left[
\left(\Gamma \vr_{n}^{\Gamma-2}+2\right)|\Grad_x \vr_{n}|^2
+\left(\Gamma b_{n}^{\Gamma-2}+2\right)|\Grad_x b_{n}|^2
\right]
\right)\dxdt \leq C,
\eeq
\begin{equation}\label{exi-23}
\int_0^{T}\int_{\Om}
\f{\ep}{\vr_{n}\vt_{n}} \f{\p p_M}{\p \vr}(\vr_{n},\vt_{n})|\Grad_x \vr_{n}|^2
\dxdt \leq C,
\eeq
\begin{equation}\label{exi-24}
C_{\ast}\vr_{n}(t,x) \leq b_{n}(t,x) \leq C^{\ast}\vr_{n}(t,x) \text{    for any    }(t,x)\in [0,T]\times \overline{\Om}.
\eeq
Here and throughout this subsection, the same letter $C$ denotes generic positive constants independent of $n$. It follows that there exists a suitable subsequence of $\{(\vr_n,b_n)\}_{n=1}^{\infty}$, not relabelled, and a weak limit $(\vr,b)$ such that
\begin{equation}\label{exi-25}
\vr_n \rightarrow \vr \text{  weakly}-\ast \text{  in }L^{\infty}(0,T;L^{\Gamma}(\Om)),
\eeq
\begin{equation}\label{exi-26}
b_n \rightarrow b \text{  weakly}-\ast \text{  in }L^{\infty}(0,T;L^{\Gamma}(\Om)),
\eeq
\begin{equation}\label{exi-27}
C_{\ast}\vr(t,x) \leq b(t,x) \leq C^{\ast}\vr(t,x) \text{    for a.e.    }(t,x)\in (0,T)\times \Om
\eeq
It is a routine matter to check, with the help of Arzela-Ascoli theorem, that$\footnote{A function $f$ belongs to $C_{weak}([0,T];L^{p}(\Om))$ if and only if it belongs to $L^{\infty}(0,T;L^{p}(\Om))$ and the mapping $t\mapsto \int_{\Om}f(t,x)\phi(x)\dx$ is continuous on $[0,T]$ for any $\phi\in  L^{p' }(\Om) $, $\f{1}{p}+\f{1}{p'}=1$.    }$
\begin{equation}\label{exi-28}
\vr_n \rightarrow \vr \text{  in }C_{weak}([0,T];L^{\Gamma}(\Om)),
\eeq
\begin{equation}\label{exi-29}
b_n \rightarrow b  \text{  in }C_{weak}([0,T];L^{\Gamma}(\Om)).
\eeq
Recalling the assumptions about viscosity coefficients (\ref{fmhd_15})$_1$, we infer from (\ref{exi-20}) and generalized Korn's inequality (\ref{ap-2}) that
\begin{equation}\label{exi-30}
\|\vu_n\|_{L^2(0,T;W_0^{1,2}(\Om))  }  \leq C      ;
\eeq
whence we may assume that
\begin{equation}\label{exi-31}
\vu_n \rightarrow \vu \text{   weakly in   }L^2(0,T;W_0^{1,2}(\Om;\R^2))
\eeq
for some $\vu \in L^2(0,T;W_0^{1,2}(\Om;\R^2))$. In view of the coercivity property of Helmholtz function (\ref{ap-3}) recalled in appendix, estimate (\ref{exi-19}), together with (\ref{fmhd_16})-(\ref{fmhd_18}), we get
\begin{equation}\label{exi-32}
\|\vt_n\|_{L^{\infty}(0,T;L^4(\Om))  }  \leq C      ;
\eeq
whence
\begin{equation}\label{exi-33}
\vt_n \rightarrow \vt \text{  weakly}-\ast \text{  in }L^{\infty}(0,T;L^{4}(\Om))
\eeq
for a suitable subsequence. Moreover, it follows directly from (\ref{exi-20}) that
\begin{equation}\label{exi-34}
\| \Grad_x \vt_n^{ \f{\Gamma}{2}  } \|_{L^2((0,T)\times\Om)}\leq C,\,\,\|  \Grad_x \vt_n^{ -\f{1}{2} }    \|_{L^2((0,T)\times\Om)}\leq C.
\eeq
Combining (\ref{exi-32}), (\ref{exi-34}) and generalized Poincar\'{e} inequality (\ref{ap-4}),
\begin{equation}\label{exi-35}
\|\vt_n\|_{L^2(0,T;W^{1,2}(\Om))}\leq C,\,\,
\|  \vt_n^{ \f{\Gamma}{2}  } \|_{L^2(0,T;W^{1,2}(\Om))}\leq C.
\eeq
Thus we may furthermore assume that
\begin{equation}\label{exi-36}
\vt_n \rightarrow \vt \text{    weakly in   } L^{2}(0,T;W^{1,2}(\Om)).
\eeq
By invoking Sobolev's embedding inequality $W^{1,2}(\Om)\hookrightarrow L^p(\Om),1\leq p< \infty$, (\ref{exi-35})$_2$ implies
\begin{equation}\label{exi-37}
\|\vt_n\|_{L^{\Gamma}(0,T;L^{p \Gamma }(\Om))  }  \leq C \,\, \text{   for any }1\leq p< \infty.
\eeq
In view of (\ref{exi-21}) and the convexity of mapping $\vt\mapsto \vt^{-3},\vt>0$, we have
\begin{equation}\label{exi-38}
\int_0^T \int_{\Om} \f{1}{\vt^3}\dxdt
\leq \liminf_{n\rightarrow \infty} \int_0^T \int_{\Om} \f{1}{\vt_n^3}\dxdt\leq C;
\eeq
whence $\vt(t,x)>0$ a.e. in $(0,T)\times \Om$ follows. We thus conclude from (\ref{exi-37}) and (\ref{exi-38}) that
\begin{equation}\label{exi-39}
\|  \log \vt_n \|_{L^p( (0,T)\times \Om )   } \leq C \,\, \text{   for any }1\leq p< \infty.
\eeq

\subsubsection{Strong convergence of $\{(\vr_n,b_n,\Grad_x \vr_n,\Grad_x b_n)\}_{n \geq 1}$}\label{f-g-1}

Due to the high nonlinearities in approximate system (\ref{exi-2})-(\ref{exi-5}), strong convergences of approximate densities and magnetic fields are indispensable. To begin with, we show strong convergence of $\{\vr_n \}_{n\geq 1}$. Notice that, in view of (\ref{fmhd_16}), (\ref{exi-19}) and (\ref{exi-37}),
\[
\left|
\int_0^T \int_{\Om} p(\vr_n,\vt_n)\Div_x \vu_n \dxdt
\right|
\leq C \int_0^T \int_{\Om} \left( \vr_n^{\gamma}+\vr_n \vt_n+\vt_n^4 \right)|\Div_x \vu_n| \dxdt
\]
\[
\leq C
\int_0^T \int_{\Om} \left( \vr_n^{\gamma}+\vr_n^{4/3} +\vt_n^4 \right)|\Grad_x \vu_n| \dxdt
\leq C;
\]
whence it follows from (\ref{exi-13}) that
\begin{equation}\label{exi-40}
\int_0^{T}\int_{\Om}
\Big\{
\ep |\Grad_x b_n|^2+
\ep \delta \left[
\left(\Gamma \vr_n^{\Gamma-2}+2\right)|\Grad_x \vr_n|^2
+\left(\Gamma b_n^{\Gamma-2}+2\right)|\Grad_x b_n|^2
\right]
\Big\} \dxdt \leq C
\eeq
for suitably large but fixed parameter $\Gamma>1$. Application of generalized Poincar\'{e} inequality (\ref{ap-4}) again yields
\begin{equation}\label{exi-41}
\|\vr_n\|_{L^2(0,T;W^{1,2}(\Om))}\leq C,\,\,
\|  \vr_n^{ \f{\Gamma}{2}  } \|_{L^2(0,T;W^{1,2}(\Om))}\leq C.
\eeq
By the same token as for $\vt_n$, it holds
\begin{equation}\label{exi-42}
\|\vr_n\|_{L^{\Gamma}(0,T;L^{p \Gamma }(\Om))  }  \leq C \,\, \text{   for any }1\leq p< \infty.
\eeq
In order to prove the strong convergence of $\{\vr_n\}_{n\geq 1}$, it suffices to obtain uniform boundedness for $\{\p_t \vr_n\}_{n\geq 1}$ and $\{\p_{ij}\vr_n\}_{n\geq 1}$ in $L^p((0,T)\times \Om),p>1$ . To do this, we first notice that (\ref{exi-30}) and (\ref{exi-40}) imply
\begin{equation}\label{exi-43}
\|\vu_n \cdot \Grad_x \vr_n \|_{L^{1}(0,T;L^{q }(\Om))  }  \leq C ,\text{     with   }\f{1}{q}=
\f{1}{2}+\f{1}{p},
\eeq
for any $2<p<\infty$. In addition, multiplying (\ref{exi-2}) by $\log \vr_n +1$ both sides gives the uniform bound
\begin{equation}\label{exi-44}
\ep \left\|\f{\Grad_x \vr_n}{\sqrt{\vr_n}}\right \|^2_{L^2((0,T)\times\Om)}\leq C.
\eeq
Consequently,
\begin{equation}\label{exi-45}
\|\vu_n \cdot \Grad_x \vr_n  \|_{L^2(0,T;L^1(\Om))}\leq
\left\|\f{\Grad_x \vr_n}{\sqrt{\vr_n}}\right \|_{L^2((0,T)\times\Om)} \|\sqrt{\vr_n}\vu_n\|_{L^{\infty}(0,T;L^2(\Om))}\leq
C.
\eeq
where we employyed the uniform boundedness of $\|\sqrt{\vr_n}\vu_n\|_{L^{\infty}(0,T;L^2(\Om))}$ due to (\ref{exi-19}). We then infer from (\ref{exi-43}), (\ref{exi-45}) and the interpolation inequality that
\begin{equation}\label{exi-46}
\|\vu_n \cdot \Grad_x \vr_n \|_{L^{ \tilde{p} }(0,T;L^{ \tilde{q} }(\Om))  }  \leq C ,\text{     for any   }  \tilde{q}\in (1,q)  , \tilde{p}=\tilde{p}(\tilde{q} )\in (1,2).
\eeq
With (\ref{exi-46}) at hand, we may invoke the maximal regularity theory for parabolic equations so as to get
\begin{equation}\label{exi-47}
\|\p_t \vr_n\| _{L^{ \tilde{p} }(0,T;L^{ \tilde{q} }(\Om))  }  \leq C ,\,\,   \|\p_{x_i x_j}\vr_n\|_{L^{ \tilde{p} }(0,T;L^{ \tilde{q} }(\Om))  }  \leq C
\eeq
for the same $\tilde{p},\tilde{q}$ as above. Recalling the compact embedding theorem, we conclude from (\ref{exi-47}) that
\begin{equation}\label{exi-48}
\vr_n \rightarrow \vr \text{     a.e. in  }(0,T)\times \Om.
\eeq
Now it is easy to pass to the limit $n\rightarrow \infty $ in (\ref{exi-2}), obtaining
\begin{equation}\label{exi-49}
\left\{\begin{aligned}
& \p_t \vr +\Div_x (\vr\vu)=\ep \Delta \vr,\,\,\text{  a.e. in }(0,T)\times \Om, \\
&  \Grad_x \vr \cdot \mathbf{n}|_{\p \Om}=0, \\
& \vr|_{t=0}=\vr_{0,\delta}. \\
\end{aligned}\right.
\end{equation}

Pointwise convergence of $\{b_n\}_{n\geq 1}$ is verified in the same strategy, with the exception that (\ref{exi-45}) must be replaced by
\[
\|\vu_n \cdot \Grad_x b_n  \|_{L^2(0,T;L^1(\Om))}\leq
\left\|\f{\Grad_x b_n}{\sqrt{b_n}}\right \|_{L^2((0,T)\times\Om)} \|\sqrt{b_n}\vu_n\|_{L^{\infty}(0,T;L^2(\Om))}
\]
\begin{equation}\label{exi-50}
\leq \sqrt{ C^{\ast}  }\left\|\f{\Grad_x b_n}{\sqrt{b_n}}\right \|_{L^2((0,T)\times\Om)}
\|\sqrt{\vr_n}\vu_n\|_{L^{\infty}(0,T;L^2(\Om))} \leq C,
\eeq
where the domination condition (\ref{exi-27}) was essentially exploited. Other calculations are left to the interested reader. Therefore, we arrive at
\begin{equation}\label{exi-51}
b_n \rightarrow b \text{     a.e. in  }(0,T)\times \Om;
\eeq
\begin{equation}\label{exi-52}
\left\{\begin{aligned}
& \p_t b +\Div_x (b \vu)=\ep \Delta b,\,\,\text{  a.e. in }(0,T)\times \Om, \\
&  \Grad_x b \cdot \mathbf{n}|_{\p \Om}=0, \\
& b|_{t=0}=b_{0,\delta}. \\
\end{aligned}\right.
\end{equation}
To proceed, we prove strong convergence of $\{\Grad_x b_n\}_{n \geq 1}$. Indeed, multiplying (\ref{exi-5}) by $2 b_n$ and integrating over $(0,\tau)\times \Om$ implies
\[
\int_{\Om} b_n^2 (\tau,x)\dx
+2 \ep \int_0^{\tau}\int_{\Om}
|\Grad_x  b_n |^2\dxdt
=\int_{\Om} b_{0,\delta}^2 \dx
-\int_0^{\tau}\int_{\Om} b_n^2 \Div_x \vu_n \dxdt,
\]
the right-hand side of which tends to
\[
\int_{\Om} b_{0,\delta}^2 \dx
-\int_0^{\tau}\int_{\Om} b^2 \Div_x \vu \dxdt
\]
according to (\ref{exi-31}), (\ref{exi-51}). Similarly, we multiply the limit equaion (\ref{exi-52})$_1$ by $2b$ to see
\[
\int_{\Om} b^2(\tau,x)\dx +2\ep \int_0^{\tau}\int_{\Om}
|\Grad_x  b|^2\dxdt=
\int_{\Om} b_{0,\delta}^2 \dx
-\int_0^{\tau}\int_{\Om} b^2 \Div_x \vu \dxdt.
\]
The above basic facts lead to
\begin{equation}\label{exi-53}
\Grad_x b_n \rightarrow   \Grad_x b \text{     strongly in    } L^{2}((0,T)\times\Om).
\eeq
Following exactly the same line, we also obtain that
\begin{equation}\label{exi-54}
\Grad_x \vr_n \rightarrow   \Grad_x \vr \text{     strongly in    } L^{2}((0,T)\times\Om).
\eeq

\subsubsection{Strong convergence of $\{\vt_n\}_{n \geq 1}$}
In order to pass to the limit in approximate entropy equation (\ref{exi-15}), strong convergence of $\{\vt_n\}_{n \geq 1}$ is necessary. Similar to the compressible Navier-Stokes-Fourier system \cite{FN1}, this is proved with the aid of Div-Curl lemma, together with suitable uniform estimates. We will thus present the sketch here for brevity. Using Gibbs' relation (\ref{fmhd_6_1}), one finds after a straightforward calculation that (\ref{exi-15}) is equivalent to
\begin{equation}\label{exi-55}
\p_t [\vr_n s(\vr_n,\vt_n)]+\Div_x [ \mathbf{R}_n^{ (1) } ]=
R_n^{ (2) }+R_n^{ (3) },
\eeq
with
\[
\mathbf{R}_n^{ (1) } :=
\vr_n s(\vr_n,\vt_n) \vu_n
-
\left[\f{\kappa(\vt_n )}{\vt_n }+\delta \left(\vt_n ^{\Gamma-1}+\f{1}{\vt_n ^2}\right)\right] \Grad_x \vt_n
\]
\[
-\ep
\left[
\vt_n s_{M}(\vr_n,\vt_n)-e_{M}(\vr_n,\vt_n)-\f{p_M(\vr_n,\vt_n)}{\vr_n}
\right]
\f{\Grad_x \vr_n}{\vt_n},
\]
\[
R_n^{ (2) }:=
\f{1}{\vt_n}\left\{
\mathbb{S}(\vt_n,\Grad_x \vu_n):\Grad_x \vu_n
+ \left[\f{\kappa(\vt_n)}{\vt_n}+\delta \left(\vt_n^{\Gamma-1}+\f{1}{\vt_n^2}\right)\right]|\Grad_x \vt_n|^2
+\delta \f{1}{\vt_n^2}
\right\}
\]
\[
+
\f{\ep}{\vt_n}|\Grad_x b_n|^2
+\f{\ep\delta }{\vt_n} \left[
\left(\Gamma \vr_n^{\Gamma-2}+2\right)|\Grad_x \vr_n|^2
+\left(\Gamma b_n^{\Gamma-2}+2\right)|\Grad_x b_n|^2
\right]
\]
\[
+\f{\ep}{\vr_n\vt_n}\f{\p p_M}{\p \vr}(\vr_n,\vt_n)|\Grad_x \vr_n|^2,
\]
\[
R_n^{ (3) }:=
-\ep \left[
e_{M}(\vr_n,\vt_n)+\vr_n \f{\p e_M}{\p \vr}(\vr_n,\vt_n)
\right]\f{\Grad_x \vr_n \cdot \Grad_x \vt_n}{\vt_n^2}-\ep \vt_n^4.
\]

Next, by the uniform estimates established so far and the structural hypotheses (\ref{fmhd_15})-(\ref{fmhd_21}), we obtain, upon choosing $\Gamma$ suitably large, the following uniform estimates:
\begin{equation}\label{exi-56}
\|(\vr_n s(\vr_n,\vt_n),\mathbf{R}_n^{ (1) })\|_{L^p((0,T)\times \Om  )} \leq C \text{   for some }p>1,
\eeq
\begin{equation}\label{exi-57}
\| (R_n^{ (2) },R_n^{ (3) }) \|_{L^1((0,T)\times \Om  )} \leq C,
\eeq
As a consequence, we set
\[
\mathbf{f}_n:= (\vr_n s(\vr_n,\vt_n),\mathbf{R}_n^{ (1) }),\,\,
\mathbf{g}_n:= (\vt_n,0,0),
\]
and apply the Div-Curl lemma recalled in Section \ref{dicu} in the time-space domain $(0,T)\times \Om$, by noticing that $L^1((0,T)\times \Om  )$ and $L^2((0,T)\times \Om  )$ are compactly embedded into $W^{-1,\tilde{p}}((0,T)\times \Om)$ for any $1<\tilde{p}<\f{3}{2}$. It thus holds that
\begin{equation}\label{exi-58}
\overline{\vr s(\vr,\vt)\vt}=\overline{\vr s(\vr,\vt)}\,\vt,
\eeq
where and throughout this paper we use overbar to denote the weak $L^1$-limit of corresponding approximation sequences. At this stage, approximation sequences are indexed by $n$. From the definition of $s(\vr,\vt)$ we know that
\begin{equation}\label{exi-59}
\vr s(\vr,\vt)=\vr s_{M}(\vr,\vt) +\f{4}{3}a \vt^3.
\eeq
To proceed, monotonicity of $\vt \mapsto \vt^3,\vt>0$ yields
\begin{equation}\label{exi-60}
\overline{\vt^4}\geq \overline{\vt^3}\,\vt.
\eeq
Moreover, it is also clear that $s_M(\vr,\vt)$ is monotonic increasing with respect to $\vt$. Thus,
\[
0\leq [s_M(\vr_n,\vt_n)-s_M(\vr_n,\vt)](\vt_n-\vt)
=s_M(\vr_n,\vt_n)(\vt_n-\vt)-s_M(\vr_n,\vt)(\vt_n-\vt);
\]
whence the strong convergence of $\{\vr_n\}_{n\geq 1}$ and weak convergence of $\{\vt_n\}_{n\geq 1}$ imply \begin{equation}\label{exi-61}
\overline{s_M(\vr,\vt)\vt}\geq \overline{s_M(\vr,\vt)} \,\vt.
\eeq
By the same token again and (\ref{exi-61}), we see
\begin{equation}\label{exi-62}
\overline{\vr s_M(\vr,\vt)\vt}=\vr\, \overline{s_M(\vr,\vt)\vt   }
\geq \vr \,\overline{s_M(\vr,\vt)} \,\vt=\overline{\vr s_M(\vr,\vt)} \, \vt.
\eeq
Putting (\ref{exi-58})-(\ref{exi-62}) together, it is easily seen that
\[
\overline{\vt^4}=\overline{\vt^3} \,\vt,
\]
which implies for suitable subsequences
\begin{equation}\label{exi-64}
\vt_n \rightarrow \vt \text{     a.e. in   } (0,T)\times \Om.
\eeq
Here we essentially employyed the strict convexity of the mapping $\vt \mapsto \vt^3,\vt>0$.

Now we are in a position to pass to the limit $n\rightarrow \infty $ in approximate entropy equation (\ref{exi-55}). Firstly, one may use the structural hypotheses and elementary inequalities to estimate (for instance, see formula (3.107) in \cite{FN1})
\[
\ep \left|
e_{M}(\vr_n,\vt_n)+\vr_n \f{\p e_M}{\p \vr}(\vr_n,\vt_n)
\right| \left| \f{\Grad_x \vr_n \cdot \Grad_x \vt_n}{\vt_n^2}\right|
\]
\begin{equation}\label{exi-65}
\leq \f{\delta}{2} \left(\vt_n^{\Gamma-2}+\f{1}{\vt_n^3}  \right)|\Grad_x \vt_n|^2
+\f{\ep \delta}{2\vt_n}  \left(\Gamma \vr_n^{\Gamma-2} +2 \right) |\Grad_x \vr_n|^2.
\eeq
In light of relation (\ref{exi-65}), we furthermore write the approximate entropy equation (\ref{exi-55}) as
\[
\p_t [\vr_n s(\vr_n,\vt_n)]+\Div_x \left\{ \vr_n s(\vr_n,\vt_n) \vu_n
-
\left[\f{\kappa(\vt_n )}{\vt_n }+\delta \left(\vt_n ^{\Gamma-1}+\f{1}{\vt_n ^2}\right)\right] \Grad_x \vt_n
\right\}
\]
\[
-\ep \Div_x \left\{
\left[
\vt_n s_{M}(\vr_n,\vt_n)-e_{M}(\vr_n,\vt_n)-\f{p_M(\vr_n,\vt_n)}{\vr_n}
\right]
\f{\Grad_x \vr_n}{\vt_n}
\right\}
\]
\[
\geq
\f{1}{\vt_n}\left\{
\mathbb{S}(\vt_n,\Grad_x \vu_n):\Grad_x \vu_n
+ \left[\f{\kappa(\vt_n)}{\vt_n}+\f{\delta}{2} \left(\vt_n^{\Gamma-1}+\f{1}{\vt_n^2}\right)\right]|\Grad_x \vt_n|^2
+\delta \f{1}{\vt_n^2}
\right\}
\]
\[
+
\f{\ep}{\vt_n}|\Grad_x b_n|^2
+\f{\ep\delta }{2\vt_n} \left[
\left(\Gamma \vr_n^{\Gamma-2}+2\right)|\Grad_x \vr_n|^2
+\left(\Gamma b_n^{\Gamma-2}+2\right)|\Grad_x b_n|^2
\right]
\]
\begin{equation}\label{exi-66}
+\f{\ep}{\vr_n\vt_n}\f{\p p_M}{\p \vr}(\vr_n,\vt_n)|\Grad_x \vr_n|^2
-\ep \vt_n^4.
\eeq
Notice that the uniform estimates (\ref{exi-25})-(\ref{exi-39}), the strong convergences of $\{\vr_n\}_{n \geq 1}$, $\{\vt_n\}_{n \geq 1}$, $\{\Grad_x \vr_n\}_{n \geq 1}$, $\{\Grad_x b_n\}_{n \geq 1}$, together with the weak convergences of $\{\Grad_x \vt_n\}_{n \geq 1}$, $\{\Grad_x \vu_n\}_{n \geq 1}$ established so far imply

\[
 \vr_n s(\vr_n,\vt_n)\rightarrow \vr s(\vr,\vt) \text{    strongly in   }L^2((0,T)\times \Om  ) ,
\]
\[
\vr_n s(\vr_n,\vt_n) \vu_n \rightarrow \vr s(\vr,\vt)\vu \text{    weakly in   }L^1((0,T)\times \Om  ) , \]
\[
\left[\f{\kappa(\vt_n )}{\vt_n }+\delta \left(\vt_n ^{\Gamma-1}+\f{1}{\vt_n ^2}\right)\right] \Grad_x \vt_n \rightarrow \left[\f{\kappa(\vt )}{\vt }+\delta \left(\vt ^{\Gamma-1}+\f{1}{\vt ^2}\right)\right] \Grad_x \vt \text{    weakly in   }L^1((0,T)\times \Om  ) ,
\]
\[
 \left[\vt_n s_{M}(\vr_n,\vt_n)-e_{M}(\vr_n,\vt_n)-\f{p_M(\vr_n,\vt_n)}{\vr_n}
\right]\f{\Grad_x \vr_n}{\vt_n}
\]
\[
\rightarrow \left[\vt s_{M}(\vr,\vt)-e_{M}(\vr,\vt)-\f{p_M(\vr,\vt)}{\vr}
\right]\f{\Grad_x \vr}{\vt}\text{    weakly in   }L^1((0,T)\times \Om  ) ,
\]
\[
 \sqrt{ \f{\mu(\vt_n)}{\vt_n}     }\left(\Grad_x \vu_n +\Grad_x^t \vu_n -\Div_x \vu_n \mathbb{I}\right)
\]
\[
\rightarrow \sqrt{ \f{\mu(\vt)}{\vt}     }\left(\Grad_x \vu +\Grad_x^t \vu -\Div_x \vu \mathbb{I}\right)\text{    weakly in   }L^2((0,T)\times \Om  ) ,
\]
\[
 \sqrt{\f{\kappa(\vt_n)}{\vt_n^2}+\f{\delta}{2} \left(\vt_n^{\Gamma-2}+\f{1}{\vt_n^3}\right)}\Grad_x \vt_n  \rightarrow  \sqrt{\f{\kappa(\vt)}{\vt^2}+\f{\delta}{2} \left(\vt^{\Gamma-2}+\f{1}{\vt^3}\right)}\Grad_x \vt \text{ weakly in }L^2((0,T)\times \Om  ) ,
\]
\[
 \f{1}{\vt_n ^{3/2}  }\rightarrow \f{1}{\vt ^{3/2}  } \text{ weakly in }L^2((0,T)\times \Om  ) ,
\]
\[
\f{\Grad_x b_n}{\sqrt{\vt_n}}  \rightarrow \f{\Grad_x b}{\sqrt{\vt}}\text{ weakly in }L^2((0,T)\times \Om ),
\]
\[
 \sqrt{  \f{\Gamma \vr_n^{\Gamma-2}+2}{\vt_n}   } \Grad_x \vr_n \rightarrow
 \sqrt{  \f{\Gamma \vr^{\Gamma-2}+2}{\vt}   } \Grad_x \vr \text{ weakly in }L^2((0,T)\times \Om  ),
\]
\[
 \sqrt{  \f{\Gamma b_n^{\Gamma-2}+2}{\vt_n}   } \Grad_x b_n \rightarrow
 \sqrt{  \f{\Gamma b^{\Gamma-2}+2}{\vt}   } \Grad_x b \text{ weakly in }L^2((0,T)\times \Om  ) ,
\]
\[
 \f{1}{\sqrt{\vr_n\vt_n}    } \sqrt{\f{\p p_M}{\p \vr}(\vr_n,\vt_n)  }\Grad_x \vr_n \rightarrow
\f{1}{\sqrt{\vr\vt}    } \sqrt{\f{\p p_M}{\p \vr}(\vr,\vt)  }\Grad_x \vr\text{ weakly in }L^2((0,T)\times \Om  ) .
\]

As a consequence, testing (\ref{exi-66}) by any $\phi \in C_c^{\infty}([0,T)\times \overline{\Om}),\phi\geq 0$ and passing to the limit $n\rightarrow \infty$ yields
\[
\int_0^T \int_{\Om} \Big(\vr s(\vr,\vt)\p_t \phi + \vr s(\vr,\vt)\vu\cdot \Grad_x \phi\Big)
\dxdt
\]
\[
-
\int_0^T \int_{\Om}
\left[\f{\kappa(\vt )}{\vt }+\delta \left(\vt ^{\Gamma-1}+\f{1}{\vt ^2}\right)\right] \Grad_x \vt
\cdot \Grad_x \phi \dxdt
\]
\[
-\ep
\int_0^T \int_{\Om}
\left[\vt s_{M}(\vr,\vt)-e_{M}(\vr,\vt)-\f{p_M(\vr,\vt)}{\vr}
\right]\f{\Grad_x \vr}{\vt}
\cdot \Grad_x \phi \dxdt
\]
\begin{equation}\label{exi-68}
+\int_0^T \int_{\Om}
\sigma_{\ep,\delta} \phi \dxdt
\leq
\int_0^T \int_{\Om} \ep \vt^4 \phi \dxdt-\int_{\Om} \vr_{0,\delta}s(\vr_{0,\delta},\vt_{0,\delta})\phi(0,\cdot) \dx,
\eeq
where
\[
\sigma_{\ep,\delta}:=
\f{1}{\vt}\left\{
\mathbb{S}(\vt,\Grad_x \vu):\Grad_x \vu
+ \left[\f{\kappa(\vt)}{\vt}+\f{\delta}{2} \left(\vt^{\Gamma-1}+\f{1}{\vt^2}\right)\right]|\Grad_x \vt|^2
+\delta \f{1}{\vt^2}
\right\}
\]
\[
+
\f{\ep}{\vt}|\Grad_x b|^2
+\f{\ep\delta }{2\vt} \left[
\left(\Gamma \vr^{\Gamma-2}+2\right)|\Grad_x \vr|^2
+\left(\Gamma b^{\Gamma-2}+2\right)|\Grad_x b|^2
\right]
+\f{\ep}{\vr\vt}\f{\p p_M}{\p \vr}(\vr,\vt)|\Grad_x \vr|^2.
\]

Finally, we pass to the limit $n\rightarrow \infty$ in approximate momentum equation (\ref{exi-3}). Clearly,
\[
\|\vr_n\vu_n\|_{L^{\infty}(0,T;L^{\f{2\Gamma}{\Gamma+1}}(\Om))   }
\leq
\|\sqrt{\vr_n}\|_{L^{\infty}(0,T;L^{2\Gamma}(\Om))   }
\|\sqrt{\vr_n}\vu_n\|_{L^{\infty}(0,T;L^{2}(\Om))   }
\leq C;
\]
whence, recalling the strong convergence of $\{\vr_n\}_{n\geq 1}$,
\[
\vr_n\vu_n
\rightarrow \vr\vu \text{  weakly}-\ast \text{  in }L^{\infty}(0,T;L^{\f{2\Gamma}{\Gamma+1}}(\Om;\R^2)).
\]
It is a routine matter to strengthen the above convergence with the help of approximate momentum equation (\ref{exi-3})
\[
\vr_n\vu_n
\rightarrow \vr\vu \text{    in   }C_{weak}([0,T];L^{\f{2\Gamma}{\Gamma+1}}(\Om;\R^2)).
\]
Since $L^q(\Om)$ is compactly embedded into $W^{-1,2}(\Om;\R^2)$ for any $1<q<\infty$, we thus obtain
\[
\vr_n\vu_n
\rightarrow \vr\vu \text{    in   }C([0,T];W^{-1,2}(\Om;\R^2)).
\]
Using (\ref{exi-31}), one infers
\begin{equation}\label{exi-69}
\vr_n\vu_n \otimes \vu_n \rightarrow \vr\vu \otimes \vu \text{    weakly in   }L^2(0,T;L^{\tilde{q}}(\Om;\R^{2\times 2})), \tilde{q}:= \f{5q}{4q+5}, 1<q<\infty.
\eeq
Based on (\ref{exi-69}) and the convergences results obtained so far, in particular the strong convergences of $\{\vr_n\}_{n\geq 1}$, $\{b_n\}_{n\geq 1}$, $\{\vt_n\}_{n\geq 1}$, $\{\Grad_x \vr_n\}_{n\geq 1}$, we now pass to the limit $n\rightarrow \infty$ in approximate momentum equation (\ref{exi-3}) to conclude
\[
\int_0^{T}\int_{\Om}
\left[ \vr\vu \cdot \p_t \vc{\phi}+\vr \vu\otimes \vu :\Grad_x \vc{\phi}+
\left( p(\vr,\vt)+\f{1}{2}b^2 +\delta \left(\vr^{\Gamma}+\vr^2+b^{\Gamma}+b^2\right)
\right)\Div_x \vc{\phi} \right] \dxdt
\]
\[
=\int_0^{T}\int_{\Om} \Big (
\mathbb{S}(\vt,\Grad_x \vu):\Grad_x \vc{\phi}+ \ep \Grad_x \vr \cdot \Grad_x \vu \cdot \vc{\phi}\Big) \dxdt
\]
\beq\label{exi-70}
-\int_{\Om} \vr_{0,\delta}\vu_{0,\delta}\cdot \vc{\phi} (0,\cdot)    \dx,
\text{     for any    }\vc{\phi} \in C_c^{\infty}([0,T)\times \Om;\R^2).
\eeq

In addition, we know from (\ref{exi-14}) the conservation of approximate total energy:
\[
\int_{\Om} \left[\f{1}{2}\vr|\vu|^2+\vr e +\f{1}{2}b^2+
\delta \left( \f{\vr^{\Gamma}  }{\Gamma-1}+\vr^2+ \f{b^{\Gamma}  }{\Gamma-1}+b^2 \right)\right](\tau,x)  \dx
\]
\[
=
\int_{\Om} \left[   \f{1}{2}\vr_{0,\delta}|\vu_{0,\delta}|^2+\vr_{0,\delta} e(\vr_{0,\delta},\vt_{0,\delta}) +\f{1}{2}b_{0,\delta}^2
+\delta \left( \f{\vr_{0,\delta}^{\Gamma}  }{\Gamma-1}+\vr_{0,\delta}^2+ \f{b_{0,\delta}^{\Gamma}  }{\Gamma-1}+b_{0,\delta}^2 \right)
\right] \dx
\]
\begin{equation}\label{exi-71}
+\int_0^{\tau}\int_{\Om}
\left(
\delta \f{1}{\vt^2}-\ep \vt^5
\right)\dxdt
\eeq
for a.e. $\tau\in (0,T)$.

\subsection{Vanishing artificial viscosity}\label{vani}
\subsubsection{Uniform-in-$\ep$ estimates}\label{unep}
From Section \ref{f-g} we know there exists $\{(\vr_{\ep},\vu_{\ep},b_{\ep},\vt_{\ep})\}_{\ep>0}$ solving the approximate equations of density (\ref{exi-49}), magnetic field (\ref{exi-51}), entropy (\ref{exi-68}), momentum (\ref{exi-70}) and total energy (\ref{exi-71}). Moreover, as a direct consequence of estimates (\ref{exi-19})-(\ref{exi-24}), the following uniform-in-$\ep$ estimates hold:
\begin{equation}\label{exi-72}
\sup_{\tau\in (0,T)}
\int_{\Om} \left[\f{1}{2}\vr_{\ep}|\vu_{\ep}|^2+\mathcal{H}_{\overline{\vt}}(\vr_{\ep},\vt_{\ep})+
\f{1}{2}b_{\ep}^2
+
\delta \left( \f{\vr_{\ep}^{\Gamma}  }{\Gamma-1}+\vr_{\ep}^2+ \f{b_{\ep}^{\Gamma}  }{\Gamma-1}+b_{\ep}^2 \right)\right](\tau,x)  \dx\leq C,
\eeq
\begin{equation}\label{exi-73}
\int_0^{T}\int_{\Om}
\f{1}{\vt_{\ep}}\left\{
\mathbb{S}(\vt_{\ep},\Grad_x \vu_{\ep}):\Grad_x \vu_{\ep}
+ \left[\f{\kappa(\vt_{\ep})}{\vt_{\ep}}+\delta \left(\vt_{\ep}^{\Gamma-1}+\f{1}{\vt_{\ep}^2}\right)\right]|\Grad_x \vt_{\ep}|^2
\right\} \dxdt \leq C,
\eeq
\begin{equation}\label{exi-74}
\int_0^{T}\int_{\Om}\left(
\delta \f{1}{\vt_{\ep}^3}+ \ep\vt_{\ep}^5+\f{\ep}{\vt_{\ep}}|\Grad_x b_{\ep}|^2\right)
\dxdt\leq C,
\eeq
\begin{equation}\label{exi-75}
\int_0^{T}\int_{\Om}\left(
\f{\ep\delta}{\vt_{\ep}} \left[
\left(\Gamma \vr_{\ep}^{\Gamma-2}+2\right)|\Grad_x \vr_{\ep}|^2
+\left(\Gamma b_{\ep}^{\Gamma-2}+2\right)|\Grad_x b_{\ep}|^2
\right]
\right)\dxdt \leq C,
\eeq
\begin{equation}\label{exi-76}
\int_0^{T}\int_{\Om}
\f{\ep}{\vr_{\ep}\vt_{\ep}} \f{\p p_M}{\p \vr}(\vr_{\ep},\vt_{\ep})|\Grad_x \vr_{\ep}|^2
\dxdt \leq C,
\eeq
\begin{equation}\label{exi-77}
C_{\ast}\vr_{\ep}(t,x) \leq b_{\ep}(t,x) \leq C^{\ast}\vr_{\ep}(t,x) \text{    for a.e.    }(t,x)\in (0,T)\times \Om.
\eeq
Here and during this subsection we denote by $C$ generic positive constants independent of $\ep$.

In analogy with Faedo-Galerkin limit, there exists a suitable subsequence of $\{(\vr_{\ep},\vu_{\ep},b_{\ep},\vt_{\ep})\}_{\ep>0}$, not relabelled, and a weak limit $(\vr,\vu,b,\vt)$ such that
\begin{equation}\label{exi-78}
\vr_{\ep} \rightarrow \vr \text{    weakly    }-\ast \text{  in }L^{\infty}(0,T;L^{\Gamma}(\Om)),
\eeq
\begin{equation}\label{exi-79}
b_{\ep} \rightarrow b \text{     weakly   }-\ast \text{  in }L^{\infty}(0,T;L^{\Gamma}(\Om)),
\eeq
\begin{equation}\label{exi-80}
\vu_{\ep} \rightarrow \vu \text{   weakly in   }L^2(0,T;W_0^{1,2}(\Om;\R^2)),
\eeq
\begin{equation}\label{exi-81}
\vt_{\ep} \rightarrow \vt \text{  weakly}-\ast \text{  in }L^{\infty}(0,T;L^{4}(\Om)),
\eeq
\begin{equation}\label{exi-82}
C_{\ast}\vr(t,x) \leq b(t,x) \leq C^{\ast}\vr(t,x) \text{    for a.e.    }(t,x)\in (0,T)\times \Om.
\eeq
Multiplying (\ref{exi-49}) by $\vr_{\ep}$ and using the uniform bounds for $\vr_{\ep},\vu_{\ep}$ gives rise to
\begin{equation}\label{exi-83}
\ep \|\Grad_x \vr_{\ep}\|_{L^2( (0,T)\times \Om  )}^2\leq C,
\eeq
which immediately implies
\begin{equation}\label{exi-84}
\ep \Grad_x \vr_{\ep} \rightarrow \mathbf{0} \text{     strongly in   } L^2( (0,T)\times \Om  ) ;
\eeq
similarly, we multipy (\ref{exi-51}) by $b_{\ep}$ to find that
\begin{equation}\label{exi-85}
\ep \|\Grad_x b_{\ep}\|_{L^2( (0,T)\times \Om  )}^2\leq C,
\eeq
\begin{equation}\label{exi-86}
\ep \Grad_x b_{\ep} \rightarrow \mathbf{0} \text{     in   } L^2( (0,T)\times \Om  ) .
\eeq

Now we pass to the limit $\ep \rightarrow 0$ in approximate equations of density and magnetic field. As in the previous section, convergences (\ref{exi-78})-(\ref{exi-79}) may be strengthened to
\begin{equation}\label{exi-87}
\vr_{\ep} \rightarrow \vr  \text{    in }C_{weak}([0,T];L^{\Gamma}(\Om)),
\eeq
\begin{equation}\label{exi-88}
b_{\ep} \rightarrow b \text{    in }C_{weak}([0,T];L^{\Gamma}(\Om)),
\eeq
with the help of Arzela-Ascoli theorem and (\ref{exi-78})-(\ref{exi-80}). Combining (\ref{exi-87})-(\ref{exi-88}) with (\ref{exi-80}),
\begin{equation}\label{exi-89}
\vr_{\ep} \vu_{\ep} \rightarrow \vr \vu \text{    weakly    }-\ast \text{  in }L^{\infty}(0,T;L^{ \f{2\Gamma}{\Gamma+1}}(\Om)),
\eeq
\begin{equation}\label{exi-90}
b_{\ep} \vu_{\ep} \rightarrow b \vu \text{    weakly    }-\ast \text{  in }L^{\infty}(0,T;L^{ \f{2\Gamma}{\Gamma+1}}(\Om)).
\eeq
In view of (\ref{exi-78}), (\ref{exi-84}) and (\ref{exi-89}), we pass to the limit $\ep\rightarrow 0$ in
(\ref{exi-49}) to conclude
\beq\label{exi-91}
\int_0^{T}\int_{\Om} \Big( \vr \p_t \phi +\vr \vu \cdot \Grad_x \phi\Big)\dxdt
+  \int_{\Om} \vr_{0,\delta}\phi(0,\cdot) \dx = 0,
\eeq
for any $\phi \in C_c^{\infty}( [0,T)\times \overline{\Om}  )$. In a similar manner,
\beq\label{exi-92}
\int_0^{T}\int_{\Om} \Big( b \p_t \phi +b \vu \cdot \Grad_x \phi \Big)\dxdt
+  \int_{\Om} b_{0,\delta}\phi(0,\cdot) \dx = 0,
\eeq
for any $\phi \in C_c^{\infty}( [0,T)\times \overline{\Om}  )$. For future analysis, it should be pointed out that the integral identities (\ref{exi-91})-(\ref{exi-92}) hold in $(0,T)\times \R^2$ provided $(\vr,\vu,b)$ are set to be zero outside $\Om$. This is a direct consequence of the estimates for $(\vr,\vu,b)$ and the non-slip boundary condition of velocity. The detailed proof can be found in Lemma 3.3 in \cite{FNP}.

\subsubsection{Strong convergence of $\{\vt_{\ep}\}_{\ep>0}$}\label{strovt}
At this level of approximation, one has to firstly show strong convergence of $\{\vt_{\ep}\}_{\ep>0}$, then of $\{\vr_{\ep}\}_{\ep>0}$. Similar to (\ref{exi-34})-(\ref{exi-39}), we notice that (\ref{exi-72})-(\ref{exi-75}) give rise to uniform-in-$\ep$ estimates:
\begin{equation}\label{exi-93}
\|\vt_{\ep}\|_{L^2(0,T;W^{1,2}(\Om))}\leq C,\,\,
\|  \vt_{\ep}^{ \f{\Gamma}{2}  } \|_{L^2(0,T;W^{1,2}(\Om))}\leq C.
\eeq
\begin{equation}\label{exi-94}
\| \vt_{\ep}^{-1} \|_{L^3((0,T)\times \Om)}\leq C,\,\,\|  \Grad_x \vt_{\ep}^{ -\f{1}{2} }    \|_{L^2((0,T)\times \Om)}\leq C,
\eeq
\begin{equation}\label{exi-95}
\|\vt_{\ep}\|_{L^{\Gamma}(0,T;L^{p \Gamma }(\Om))  }  \leq C \,\, \text{   for any }1\leq p< \infty,
\eeq
\begin{equation}\label{exi-96}
\|  \log \vt_{\ep} \|_{L^p( (0,T)\times \Om )   } \leq C \,\, \text{   for any }1\leq p< \infty,
\eeq
\begin{equation}\label{exi-97}
\|\sqrt{\ep}\vr_{\ep}\|_{L^2(0,T;W^{1,2}(\Om))}\leq C,\,\,
\| \sqrt{\ep} \vr_{\ep}^{ \f{\Gamma}{2}  } \|_{L^2(0,T;W^{1,2}(\Om))}\leq C.
\eeq
In order to apply the Div-Curl lemma conveniently, we reformulate the approximate entropy equation by introducing non-negative regular Borel measure$\footnote{ Let $Q$ be a locally compact Hausdorff space. In accordance with Riesz representation theorem, a non-negative regular Borel measure can be identified with a non-negative linear functional on $C_c(Q)$. All such measures constitute a vector space, denoted by $\mathcal{M}^{+}(Q)$. }$:
\[
\left< \Sigma_{\ep,\delta};\phi \right>=
\int_0^T \int_{\Om} \ep \vt_{\ep}^4 \phi \dxdt-\int_{\Om} \vr_{0,\delta}s(\vr_{0,\delta},\vt_{0,\delta})\phi(0,\cdot) \dx
\]
\[
-\int_0^T \int_{\Om} \Big(\vr_{\ep} s(\vr_{\ep},\vt_{\ep})\p_t \phi + \vr_{\ep} s(\vr_{\ep},\vt_{\ep})\vu_{\ep}\cdot \Grad_x \phi\Big)
\dxdt
\]
\[
+
\int_0^T \int_{\Om}
\left[\f{\kappa(\vt_{\ep} )}{\vt_{\ep} }+\delta \left(\vt_{\ep} ^{\Gamma-1}+\f{1}{\vt_{\ep} ^2}\right)\right] \Grad_x \vt_{\ep}
\cdot \Grad_x \phi \dxdt
\]
\begin{equation}\label{exi-98}
+\ep
\int_0^T \int_{\Om}
\left[\vt _{\ep} s_{M}(\vr_{\ep},\vt_{\ep})-e_{M}(\vr_{\ep},\vt_{\ep})-\f{p_M(\vr_{\ep},\vt_{\ep})}{\vr_{\ep}}
\right]\f{\Grad_x \vr_{\ep}}{\vt_{\ep}}
\cdot \Grad_x \phi \dxdt,
\eeq
for any $\phi \in C_c^{\infty}([0,T)\times \overline{\Om}),\phi\geq 0$. Alternatively, we may write (\ref{exi-98}) in the weak form
\begin{equation}\label{exi-99}
\p_t [\vr_{\ep} s(\vr_{\ep},\vt_{\ep})]+\Div_x [ \mathbf{Z}_{\ep}^{ (1) } ]=
Z_n^{ (2) }+Z_{\ep}^{ (3) },
\eeq
with

\[
\mathbf{Z}_{\ep}^{ (1) }:=
 \vr_{\ep} s(\vr_{\ep},\vt_{\ep}) \vu_{\ep}
-
\left[\f{\kappa(\vt_{\ep} )}{\vt_{\ep} }+\delta \left(\vt_{\ep} ^{\Gamma-1}+\f{1}{\vt_{\ep}^2}\right)\right] \Grad_x \vt_{\ep}
\]
\[
-\ep
\left[
\vt_{\ep} s_{M}(\vr_{\ep},\vt_{\ep})-e_{M}(\vr_{\ep},\vt_{\ep})-\f{p_M(\vr_{\ep},\vt_{\ep})}{\vr_{\ep}}
\right]
\f{\Grad_x \vr_{\ep}}{\vt_{\ep}},
\]
\[
Z_{\ep}^{ (2) }:=\Sigma_{\ep,\delta},\,\,Z_{\ep}^{ (3) }:=-\ep \vt_{\ep}^4.
\]

Now we are in a position to apply the Div-Curl lemma recalled in Section \ref{dicu} in the time-space domain $(0,T)\times \Om$ with
\[
\mathbf{f}_{\ep}:= (\vr_{\ep} s(\vr_{\ep},\vt_{\ep}),\mathbf{Z}_{\ep}^{ (1) }),\,\,
\mathbf{g}_{\ep}:= (F(\vt_{\ep}),0,0), \,\,  F \in W^{1,\infty}(\R^{+}).
\]
Using the uniform-in-$\ep$ estimates obtained, we know that
\begin{equation}\label{exi-100}
Z_{\ep}^{ (3) } \rightarrow 0 \text{   strongly in   }L^1((0,T)\times \Om  ).
\eeq
Seeing next that $L^1((0,T)\times \Om  )$ and $\mathcal{M}^{+}( [0,T]\times \overline{\Om}  )$ are compactly embedded into $W^{-1,\tilde{p}}((0,T)\times \Om)$ for any $1<\tilde{p}<\f{3}{2}$. Similar as before,
\begin{equation}\label{exi-101}
\|(\vr_{\ep} s(\vr_{\ep},\vt_{\ep}),\mathbf{Z}_{\ep}^{ (1) })\|_{L^p((0,T)\times \Om  )} \leq C \text{   for some }p>1,
\eeq
and moreover,
\begin{equation}\label{exi-102}
\ep
\left[
\vt_{\ep} s_{M}(\vr_{\ep},\vt_{\ep})-e_{M}(\vr_{\ep},\vt_{\ep})-\f{p_M(\vr_{\ep},\vt_{\ep})}{\vr_{\ep}}
\right]
\f{\Grad_x \vr_{\ep}}{\vt_{\ep}}
\rightarrow 0 \text{  strongly in }L^p((0,T)\times \Om  ),
\eeq
for some $p>1$. Consequently,
\begin{equation}\label{exi-103}
\overline{\vr s(\vr,\vt) F(\vt)}=\overline{\vr s(\vr,\vt)}\,\, \overline{ F(\vt)},
\eeq
with $F\in W^{1,\infty}(\R^+)$. In the same spirit of Faedo-Galerkin limit, it is shown that$\footnote{The proof is based on the tool of parameterized Young measures and we refer to Section 3.6.2 in \cite{FN1} for similar calculations. It should be remarked that strong convergence of $\{\vr_{\ep}\}_{\ep>0}$ has not been proved yet, making the proof of these inequalities much more delicate. }$:
\begin{equation}\label{exi-104}
\left\{\begin{aligned}
& \overline{  \vr s_M(\vr,\vt)F(\vt)  } \geq \overline{ \vr s_M(\vr,\vt)    }\,\overline{F(\vt)   }, \\
& \overline{  \vt^3 F(\vt)  } \geq \overline{\vt^3} \,\overline{F(\vt)   },   \\
\end{aligned}\right.
\end{equation}
for any monotonic increasing $F\in W^{1,\infty}(\R^+)$. It follows from (\ref{exi-103})-(\ref{exi-104}) that
\begin{equation}\label{exi-105}
\vt_{\ep} \rightarrow \vt \text{     a.e. in   } (0,T)\times \Om
\eeq
for suitable subsequences. In addition, similar to (\ref{exi-38}),
\begin{equation}\label{exi-106}
\int_0^T \int_{\Om} \f{1}{\vt^3}\dxdt
\leq \liminf_{n\rightarrow \infty} \int_0^T \int_{\Om} \f{1}{\vt_{\ep}^3}\dxdt\leq C;
\eeq
whence $\vt(t,x)>0$ a.e. in $(0,T)\times \Om$.

Now we pass to the limit $\ep\rightarrow 0$ in the approximate entropy equation (\ref{exi-98}). Application of Div-Curl lemma again shows
\begin{equation}\label{exi-107}
\vr_{\ep} s(\vr_{\ep},\vt_{\ep}) \vu_{\ep}
\rightarrow   \overline{\vr s(\vr,\vt)}  \,\vu \text{   weakly in   }L^p((0,T)\times \Om  )
\eeq
for some $p>1$. Clearly, there exists $\Sigma^{(\delta)}\in \mathcal{M}^{+}( [0,T]\times \overline{\Om}  ) $ such that
\begin{equation}\label{exi-108}
\Sigma_{\ep,\delta} \rightarrow  \Sigma^{(\delta)} \text{   weakly}-\ast \text{ in }\mathcal{M}( [0,T]\times \overline{\Om}  )
\eeq
as $\ep\rightarrow 0$. Consequently, with the help of (\ref{exi-81}), (\ref{exi-100}), (\ref{exi-102}), (\ref{exi-105}), (\ref{exi-107}) and (\ref{exi-108}), we let $\ep\rightarrow 0$ in (\ref{exi-98}) to infer that
\[
\int_0^T \int_{\Om}  \Big(\overline{\vr s(\vr,\vt)}\p_t \phi
+\overline{\vr s(\vr,\vt)}\,\vu \cdot \Grad_x \phi \Big)\dxdt
\]
\[
-
\int_0^T \int_{\Om}
\left[\f{\kappa(\vt )}{\vt }+\delta \left(\vt ^{\Gamma-1}+\f{1}{\vt ^2}\right)\right] \Grad_x \vt
\cdot \Grad_x \phi \dxdt
\]
\begin{equation}\label{exi-109}
+\left<\Sigma^{(\delta)} ;\phi\right> = -\int_{\Om} \vr_{0,\delta}s(\vr_{0,\delta},\vt_{0,\delta})\phi(0,\cdot) \dx
\eeq
for any $\phi \in C_c^{\infty}([0,T)\times \overline{\Om}),\phi\geq 0$. Furthermore, observe that all $\ep$-dependent terms in $\sigma_{\ep,\delta}$ are non-negative, yielding
\begin{equation}\label{exi-110}
\Sigma^{(\delta)}\geq
\f{1}{\vt}\left\{
\mathbb{S}(\vt,\Grad_x \vu):\Grad_x \vu
+ \left[\f{\kappa(\vt)}{\vt}+\f{\delta}{2} \left(\vt^{\Gamma-1}+\f{1}{\vt^2}\right)\right]|\Grad_x \vt|^2
+\delta \f{1}{\vt^2}
\right\},
\eeq
to be understood in $ \mathcal{M}^{+}( [0,T]\times \overline{\Om}  ) $.

\subsubsection{Strong convergence of $\{\vr_{\ep}\}_{\ep>0}$}\label{strovr}
In order to pass to the limit $\ep\rightarrow 0$ in the approximate momentum equation (\ref{exi-70}), strong convergences of $\{\vr_{\ep}\}_{\ep>0}$ and $\{b_{\ep}\}_{\ep>0}$ are indispensable due to the nonlinearities in pressure. Before this, it is also necessary to improve their integrability uniformly in $\ep$. Indeed, similar to compressible Navier-Stokes system \cite{LS} and Navier-Stokes-Fourier system \cite{FN1}, this step is finished with the help of Bogovskii operator (see \cite{BG,FN1}). More specifically, by choosing
\[
\vc{\phi}(t,x)=\psi(t)\mathcal{B}\left( \vr_{\ep}-\f{1}{|\Om|} \int_{\Om}\vr_{\ep} \dx  \right),\,\,
\psi(t)\in C_c^{\infty}( (0,T) )
\]
as a test function in the approximate momentum equation (\ref{exi-70}) and taking advantage of the uniform-in-$\ep$ estimates obtained, we arrive at
\begin{equation}\label{exi-111}
\int_0^T \psi \int_{\Om}\left( p(\vr_{\ep},\vt_{\ep})+\f{1}{2}b_{\ep}^2 +\delta \left(\vr_{\ep}^{\Gamma}+\vr_{\ep}^2+b_{\ep}^{\Gamma}+b_{\ep}^2\right)
\right) \vr_{\ep} \dxdt \leq C.
\eeq
Since the proof is similar to that of Navier-Stokes system, the details are thus omitted here. In particular, we conclude from the domination condition (\ref{exi-77}) and the structural hypothesis of pressure that
\begin{equation}\label{exi-112}
\|(\vr_{\ep},b_{\ep})\|_{L^{\Gamma+1}((0,T)\times\Om )    }\leq C,
\eeq
\begin{equation}\label{exi-113}
\|p_M(\vr_{\ep},\vt_{\ep})\|_{L^{p}((0,T)\times\Om )    }\leq C \,\, \text{      for some   }p>1.
\eeq
It follows from (\ref{exi-80}), (\ref{exi-84}), (\ref{exi-87}), (\ref{exi-105}), (\ref{exi-112}), (\ref{exi-113}) that
\begin{equation}\label{exi-114}
\left\{\begin{aligned}
& \ep \Grad_x \vr_{\ep}\cdot \Grad_x \vu_{\ep}\rightarrow 0 \text{     strongly in   } L^1( (0,T)\times \Om  ), \\
&  \vr_{\ep}\vu_{\ep} \otimes \vu_{\ep} \rightarrow \vr\vu \otimes \vu \text{    weakly in   }L^p((0,T)\times \Om), \text{  for some }p>1, \\
& \mathbb{S}(\vt_{\ep},\Grad_x \vu_{\ep}) \rightarrow
\mathbb{S}(\vt,\Grad_x \vu)\text{    weakly in   }L^p((0,T)\times \Om), \text{  for some }p>1,   \\
\end{aligned}\right.
\end{equation}
Now we pass to the limit $\ep\rightarrow 0$ in the approximate momentum equation (\ref{exi-70}), in view of (\ref{exi-114}) and (\ref{exi-105}), to conclude that
\[
\int_0^{T}\int_{\Om}
\left[  \vr\vu \cdot \p_t \vc{\phi}+\vr \vu\otimes \vu :\Grad_x \vc{\phi}+
\left( \overline{p_{M,\delta}(\vr,\vt,b)}+\f{a}{3}\vt^4 \right)  \Div_x \vc{\phi} \right] \dxdt
\]
\beq\label{exi-115}
=\int_0^{T}\int_{\Om}
\mathbb{S}(\vt,\Grad_x \vu):\Grad_x \vc{\phi} \dxdt
-\int_{\Om} \vr_{0,\delta}\vu_{0,\delta}\cdot \vc{\phi} (0,\cdot)    \dx,
\eeq
for any $\vc{\phi} \in C_c^{\infty}([0,T)\times \Om;\R^2)$. Here,
\[
p_{M,\delta}(\vr,\vt,b):=p_M(\vr,\vt)+\f{1}{2}b^2 +\delta \left(\vr^{\Gamma}+\vr^2+b^{\Gamma}+b^2\right)
\]
\[
=\vr^{\gamma}+\vr\vt+\f{1}{2}b^2 +\delta \left(\vr^{\Gamma}+\vr^2+b^{\Gamma}+b^2\right).
\]

Upon setting
\[
\zeta_{\ep}:=
\begin{cases}
\f{b_{\ep}}{\vr_{\ep}}& \text{ if } \vr_{\ep}>0, \\
\f{C_{\ast}+C^{\ast}   }{2}& \text{ if }\vr_{\ep}=0,\\
\end{cases}
\,\,\,\,
\zeta:=
\begin{cases}
\f{b}{\vr}& \text{ if } \vr>0, \\
\f{C_{\ast}+C^{\ast}   }{2}& \text{ if }\vr=0,\\
\end{cases}
\]
one may rewrite $p_{M,\delta}(\vr,\vt,b)$ in another form
\[
\mathcal{P}_{M,\delta}(\vr,\vt,\zeta)=\vr^{\gamma}+\vr\vt+\f{1}{2}\zeta^2\vr^2+
\delta \left(\vr^{\Gamma}+\vr^2+\zeta^{\Gamma}\vr^{\Gamma}+\zeta^2 \vr^2\right).
\]
Therefore,
\[
p_{M,\delta}(\vr_{\ep},\vt_{\ep},b_{\ep})
=\mathcal{P}_{M,\delta}(\vr_{\ep},\vt_{\ep},\zeta_{\ep})
\]
\beq\label{exi-116}
=\mathcal{P}_{M,\delta}(\vr_{\ep},\vt_{\ep},\zeta_{\ep})
-\mathcal{P}_{M,\delta}(\vr_{\ep},\vt,\zeta)+\mathcal{P}_{M,\delta}(\vr_{\ep},\vt,\zeta).
\eeq

We are ready to invoke the following crucial lemma from compressible two-fluid model due to firstly Vasseur et al. \cite{VWY} and then refined by Novotn\'{y} et al. \cite{NP}, indicating some kind of compactness for $\{\zeta_{\ep}\}_{\ep>0}$. With the help of this lemma, strong convergence of $\{b_{\ep}\}_{\ep>0}$ will be reduced to that of $\{\vr_{\ep}\}_{\ep>0}$ to some extent.
\begin{Lemma}\label{comp-lem}
Let the uniform-in-$\ep$ estimates (\ref{exi-77})-(\ref{exi-80}) be in force for $\{(\vr_{\ep},b_{\ep},\vu_{\ep})\}_{\ep>0}$. Then, for any $1\leq p<\infty$,
\beq\label{exi-117}
\int_{\Om}\vr_{\ep} |\zeta_{\ep}-\zeta|^p (\tau,\cdot)\dx \rightarrow 0 \text{   for any  } \tau\in [0,T].
\eeq
\end{Lemma}
The detailed proof of Lemma \ref{comp-lem} can be found in Proposition 11 \cite{NP}, which is a generalization of Lemma 2.1 \cite{VWY}.

As a direct consequence of (\ref{exi-117}), (\ref{exi-77})-(\ref{exi-79}), strong convergence of $\{\vt_{\ep}\}_{\ep>0}$, H\"{o}lder's inequality,
\[
\int_0^T\int_{\Om}\left|
\mathcal{P}_{M,\delta}(\vr_{\ep},\vt_{\ep},\zeta_{\ep})
-\mathcal{P}_{M,\delta}(\vr_{\ep},\vt,\zeta)\right|
\dxdt
\]
\[
\leq C \left(
\int_0^T\int_{\Om} \vr_{\ep}|\vt_{\ep}-\vt| \dxdt
+\int_0^T\int_{\Om} \left[
\vr_{\ep}^2+\delta \left(\vr_{\ep}^{\Gamma} +\vr_{\ep}^2   \right)
\right] |\zeta_{\ep}-\zeta| \dxdt
\right)
\]
\beq\label{exi-118}
 \rightarrow 0 \,\,\text{     as         }\ep \rightarrow 0,
\eeq
where the two-sided bounds $C_{\ast} \leq \zeta_{\ep}\leq C^{\ast}$ were also used basically. We thus infer from (\ref{exi-116}) and (\ref{exi-118}) that
\beq\label{exi-119}
\overline{p_{M,\delta}(\vr,\vt,b)}=
\overline{\overline{\mathcal{P}_{M,\delta}(\vr,\vt,\zeta)}}.
\eeq
Here, we adopted the symbol from \cite{NP} denoting by $\overline{\overline{\mathcal{P}_{M,\delta}(\vr,\vt,\zeta)}}$ the weak $L^1( (0,T)\times \Om )$ limit of $\{\mathcal{P}_{M,\delta}(\vr_{\ep},\vt,\zeta)\}_{\ep>0}$. Therefore, the limit momentum equation (\ref{exi-115}) is equivalent as
\[
\int_0^{T}\int_{\Om}
\left[  \vr\vu \cdot \p_t \vc{\phi}+\vr \vu\otimes \vu :\Grad_x \vc{\phi}+
\left( \overline{\overline{\mathcal{P}_{M,\delta}(\vr,\vt,\zeta)}}+\f{a}{3}\vt^4 \right)  \Div_x \vc{\phi} \right] \dxdt
\]
\beq\label{exi-120}
=\int_0^{T}\int_{\Om}
\mathbb{S}(\vt,\Grad_x \vu):\Grad_x \vc{\phi} \dxdt
-\int_{\Om} \vr_{0,\delta}\vu_{0,\delta}\cdot \vc{\phi} (0,\cdot)    \dx,
\eeq
for any $\vc{\phi} \in C_c^{\infty}([0,T)\times \Om;\R^2)$. To proceed, we follow the ideas from compressible Navier-Stokes system \cite{FNP,HD1,LS,SD1} by establishing the remarkable ``effective viscous flux identity". Similar arguments also appeared in compressible two-fluid models \cite{NP,VWY}.
\begin{Lemma}\label{effe}
It holds that
\[
\lim_{\ep\rightarrow 0}\int_0^T \psi \int_{\Om} \varphi
\left[
\mathcal{P}_{M,\delta}(\vr_{\ep},\vt,\zeta)-\mu(\vt_{\ep}) \Div_x \vu_{\ep}
\right]
\vr_{\ep}\dxdt
\]
\beq\label{exi-121}
=\int_0^T \psi \int_{\Om} \varphi
\left[
\overline{\overline{\mathcal{P}_{M,\delta}(\vr,\vt,\zeta) }} -\mu(\vt) \Div_x \vu
\right]
\vr \dxdt
\eeq
for any $\psi \in C_c^{\infty}((0,T)),\varphi\in C_c^{\infty}(\Om)$.
\end{Lemma}
{\bf{Proof. }} Let $\Delta^{-1}$ be the inverse of Laplacian on $\R^2$ with Fourier symbol $\f{-1}{|\xi|^2}$; the operator $\mathcal{A}_j$ is defined as $\p _{x_j}\Delta^{-1}$ with Fourier symbol $\f{-i \xi_j}{|\xi|^2}$; the operator $\mathcal{R}_{i,j}$ is defined as $\p _{x_j}\mathcal{A}_i$ with Fourier symbol $\f{\xi_i \xi_j}{|\xi|^2}$. For brevity, we set $\mathcal{A}:=(\mathcal{A}_1,\mathcal{A}_2)$ and $\mathcal{R}=\{\mathcal{R}_{i,j}\}_{i,j=1}^2$. Notice that the approximate continuity equation (\ref{exi-49})$_1$ holds in $(0,T)\times \R^2$ upon extending $(\vr_{\ep},\vu_{\ep})$ to be zero outside $\Om$, i.e.,
\beq\label{exi-122}
\p_t (1_{\Om} \vr_{\ep}) =\Div_x ( \ep  1_{\Om} \Grad_x \vr_{\ep} - 1_{\Om} \vr_{\ep}\vu_{\ep})
\text{   a.e. in  } (0,T)\times \R^2.
\eeq

Testing the approximate momentum equation (\ref{exi-70}) by $\psi\varphi \mathcal{A}(1_{\Om} \vr_{\ep})$ and making use of (\ref{exi-122}) gives rise to
\[
\int_0^T\int_{\Om} \psi \varphi
\Big(
\mathcal{P}_{M,\delta}(\vr_{\ep},\vt,\zeta)\vr_{\ep} -\mathbb{S}(\vt_{\ep},\Grad_x \vu_{\ep}):\mathcal{R}(1_{\Om} \vr_{\ep})
\Big) \dxdt
\]
\[
=- \int_0^T\int_{\Om} \psi \varphi \Big(
\mathcal{P}_{M,\delta}(\vr_{\ep},\vt_{\ep},\zeta_{\ep})
-\mathcal{P}_{M,\delta}(\vr_{\ep},\vt,\zeta)\Big)\vr_{\ep}\dxdt
- \int_0^T\int_{\Om} \psi \varphi \f{a}{3}\vt_{\ep}^4 \vr_{\ep} \dxdt
\]
\[
+
\int_0^T\int_{\Om} \psi \varphi
\Big(
\vr_{\ep}\vu_{\ep}\cdot \mathcal{R}(1_{\Om} \vr_{\ep}\vu_{\ep})-
(\vr_{\ep}\vu_{\ep}\otimes \vu_{\ep} ):\mathcal{R} (1_{\Om} \vr_{\ep})
\Big)\dxdt
\]
\[
-\ep
\int_0^T\int_{\Om} \psi \varphi
\vr_{\ep}\vu_{\ep} \cdot \mathcal{A} ( \Div_x ( 1_{\Om} \Grad_x \vr_{\ep})  )
\dxdt
\]
\[
-\int_0^T\int_{\Om} \psi
\mathcal{P}_{M,\delta}(\vr_{\ep},\vt_{\ep},\zeta_{\ep})\Grad_x \varphi \cdot \mathcal{A} ( 1_{\Om} \vr_{\ep}  )
\dxdt
\]
\[
+
\int_0^T\int_{\Om} \psi
\mathbb{S}(\vt_{\ep},\Grad_x \vu_{\ep}): [\Grad_x \varphi  \otimes \mathcal{A} ( 1_{\Om} \vr_{\ep}  )]
\dxdt
\]
\[
-
\int_0^T\int_{\Om} \psi
(\vr_{\ep}\vu_{\ep}\otimes \vu_{\ep} ): [\Grad_x \varphi  \otimes \mathcal{A} ( 1_{\Om} \vr_{\ep}  )]
\dxdt
\]
\[
-
\int_0^T\int_{\Om} \p_t \psi\varphi
\vr_{\ep}\vu_{\ep} \cdot \mathcal{A} ( 1_{\Om} \vr_{\ep}  ) \dxdt
\]
\beq\label{exi-123}
+\ep
\int_0^T\int_{\Om}
 \psi \varphi (\Grad_x \vr_{\ep}\cdot \Grad_x \vu_{\ep} ) \cdot \mathcal{A} ( 1_{\Om} \vr_{\ep}  ) \dxdt
 =: \sum_{j=1}^9 \mathcal{I}_j^{(1),\ep  }  .
\eeq

In a similar manner, we use $\psi\varphi \mathcal{A}(1_{\Om} \vr)$ as a test function in the limit momentum equation (\ref{exi-120}) to see that
\[
\int_0^T\int_{\Om} \psi \varphi \Big(
\overline{\overline{\mathcal{P}_{M,\delta}(\vr,\vt,\zeta)}}\, \vr
-\mathbb{S}(\vt,\Grad_x\vu):\mathcal{R}(1_{\Om}\vr)
\Big)\dxdt
\]
\[
=-\int_0^T\int_{\Om} \psi \varphi \f{a}{3}\vt^4\vr \dxdt
\]
\[
+
\int_0^T\int_{\Om} \psi \varphi
\Big(
\vr\vu\cdot \mathcal{R}(1_{\Om} \vr\vu)-
(\vr\vu\otimes \vu ):\mathcal{R} (1_{\Om} \vr)
\Big)\dxdt
\]
\[
-
\int_0^T\int_{\Om} \psi \,
\overline{\overline{\mathcal{P}_{M,\delta}(\vr,\vt,\zeta)}}\,
\Grad_x \varphi \cdot \mathcal{A} ( 1_{\Om} \vr  )
\dxdt
\]
\[
+
\int_0^T\int_{\Om} \psi
\mathbb{S}(\vt,\Grad_x \vu): [\Grad_x \varphi  \otimes \mathcal{A} ( 1_{\Om} \vr  )]
\dxdt
\]
\[
-
\int_0^T\int_{\Om} \psi
(\vr\vu\otimes \vu ): [\Grad_x \varphi  \otimes \mathcal{A} ( 1_{\Om} \vr  )]
\dxdt
\]
\beq\label{exi-124}
-
\int_0^T\int_{\Om} \p_t \psi\varphi
\vr\vu \cdot \mathcal{A} ( 1_{\Om} \vr  ) \dxdt
 =: \sum_{j=1}^6 \mathcal{I}_j^{(2)  }  .
\eeq

Now we analyze the right-hand side of (\ref{exi-123}) term by term in the regime of vanishing artificial viscosity. Based on (\ref{exi-118}) and (\ref{exi-78}), it holds that $\mathcal{I}_1^{(1),\ep  }\rightarrow 0$ as $\ep \rightarrow 0$. It follows from H\"{o}rmander-Mikhlin theorem, Sobolev's embedding inequality and (\ref{exi-87}) that
\beq\label{exi-125}
\mathcal{A} ( 1_{\Om} \vr_{\ep}  )
\rightarrow
\mathcal{A} ( 1_{\Om} \vr  ) \text{    in   } C([0,T]\times \overline{\Om};\R^2).
\eeq
From (\ref{exi-116}), (\ref{exi-118}) and (\ref{exi-125}) we see $\mathcal{I}_5^{(1),\ep  }\rightarrow \mathcal{I}_3^{(2)  } $ as $\ep \rightarrow 0$. Strong convergence of $\{\vt_{\ep}\}_{\ep>0}$ (\ref{exi-105}) and weak convergence of $\{\vr_{\ep}\}_{\ep>0}$ (\ref{exi-78}) together imply that $\mathcal{I}_2^{(1),\ep  }\rightarrow \mathcal{I}_1^{(2)  } $ as $\ep \rightarrow 0$. We deduce easily from (\ref{exi-80}), (\ref{exi-84}) and (\ref{exi-125}) that $\mathcal{I}_9^{(1),\ep  }\rightarrow 0$ as $\ep \rightarrow 0$. Similarly, $\mathcal{I}_4^{(1),\ep  }\rightarrow 0$ as $\ep \rightarrow 0$. Following the same line, we essentially use the convergences (\ref{exi-78})-(\ref{exi-81}), (\ref{exi-105}) and (\ref{exi-125}) in order to conclude that $\mathcal{I}_6^{(1),\ep  }\rightarrow \mathcal{I}_4^{(2)  } $, $\mathcal{I}_7^{(1),\ep  }\rightarrow \mathcal{I}_5^{(2)  } $, $\mathcal{I}_8^{(1),\ep  }\rightarrow \mathcal{I}_6^{(2)  } $ as $\ep \rightarrow 0$. Consequently, upon passing to the limit $\ep\rightarrow 0$ in (\ref{exi-123}) and using the convergence results obtained, (\ref{exi-124}),
\[
\lim_{\ep\rightarrow 0}
 \int_0^T\int_{\Om} \psi \varphi
\Big(
\mathcal{P}_{M,\delta}(\vr_{\ep},\vt,\zeta)\vr_{\ep} -\mathbb{S}(\vt_{\ep},\Grad_x \vu_{\ep}):\mathcal{R}(1_{\Om} \vr_{\ep})
\Big) \dxdt
\]
\[
=
\int_0^T\int_{\Om} \psi \varphi \Big(
\overline{\overline{\mathcal{P}_{M,\delta}(\vr,\vt,\zeta)}}\, \vr
-\mathbb{S}(\vt,\Grad_x\vu):\mathcal{R}(1_{\Om}\vr)
\Big)\dxdt
\]
\[
+\lim_{\ep\rightarrow 0}
\int_0^T\int_{\Om} \psi \varphi
\Big(
\vr_{\ep}\vu_{\ep}\cdot \mathcal{R}(1_{\Om} \vr_{\ep}\vu_{\ep})-
(\vr_{\ep}\vu_{\ep}\otimes \vu_{\ep} ):\mathcal{R} (1_{\Om} \vr_{\ep})
\Big)\dxdt
\]
\beq\label{exi-126}
-
\int_0^T\int_{\Om} \psi \varphi
\Big(
\vr\vu\cdot \mathcal{R}(1_{\Om} \vr\vu)-
(\vr\vu\otimes \vu ):\mathcal{R} (1_{\Om} \vr)
\Big)\dxdt.
\eeq
In analogy with compressible Navier-Stokes system \cite{FNP}, see Lemma 3.4 therein, one may utilize another version of Div-Curl lemma so as to conclude that (the details are thus omitted):
\[
\lim_{\ep\rightarrow 0}
\int_0^T\int_{\Om} \psi \varphi
\Big(
\vr_{\ep}\vu_{\ep}\cdot \mathcal{R}(1_{\Om} \vr_{\ep}\vu_{\ep})-
(\vr_{\ep}\vu_{\ep}\otimes \vu_{\ep} ):\mathcal{R} (1_{\Om} \vr_{\ep})
\Big)\dxdt
\]
\[
=
\int_0^T\int_{\Om} \psi \varphi
\Big(
\vr\vu\cdot \mathcal{R}(1_{\Om} \vr\vu)-
(\vr\vu\otimes \vu ):\mathcal{R} (1_{\Om} \vr)
\Big)\dxdt;
\]
whence (\ref{exi-126}) reduces to
\[
\lim_{\ep\rightarrow 0}
 \int_0^T\int_{\Om} \psi \varphi
\Big(
\mathcal{P}_{M,\delta}(\vr_{\ep},\vt,\zeta)\vr_{\ep} -\mathbb{S}(\vt_{\ep},\Grad_x \vu_{\ep}):\mathcal{R}(1_{\Om} \vr_{\ep})
\Big) \dxdt
\]
\beq\label{exi-127}
=
\int_0^T\int_{\Om} \psi \varphi \Big(
\overline{\overline{\mathcal{P}_{M,\delta}(\vr,\vt,\zeta)}}\, \vr
-\mathbb{S}(\vt,\Grad_x\vu):\mathcal{R}(1_{\Om}\vr)
\Big)\dxdt.
\eeq
To proceed, notice that
\[
\int_0^T\int_{\Om} \psi \varphi
\mathbb{S}(\vt_{\ep},\Grad_x \vu_{\ep}):\mathcal{R}(1_{\Om} \vr_{\ep})
\dxdt
\]
\[
=
\int_0^T\int_{\Om} \psi \varphi
\mu(\vt_{\ep})(\Grad_x \vu_{\ep}+\Grad_x^t \vu_{\ep}-\Div_x \vu_{\ep}\mathbb{I} ): \mathcal{R}(1_{\Om} \vr_{\ep})
\dxdt
\]
\[=
\int_0^T\int_{\Om} \psi \mathcal{R}:\left[
\varphi\mu(\vt_{\ep})(\Grad_x \vu_{\ep}+\Grad_x^t \vu_{\ep})
\right] \vr_{\ep}
\dxdt
-
\int_0^T\int_{\Om} \psi \varphi
\mu(\vt_{\ep}) \vr_{\ep}\Div_x \vu_{\ep}\dxdt
\]
\beq\label{exi-128}
=
\int_0^T\int_{\Om} \psi \varphi
\mu(\vt_{\ep}) \vr_{\ep}\Div_x \vu_{\ep}\dxdt
+\int_0^T\int_{\Om} \psi \vr_{\ep} \omega(\vt_{\ep},\vu_{\ep})\dxdt,
\eeq
where
\[
\omega(\vt_{\ep},\vu_{\ep}):=
\mathcal{R}:\left[
\varphi\mu(\vt_{\ep})(\Grad_x \vu_{\ep}+\Grad_x^t \vu_{\ep})
\right]-\varphi\mu(\vt_{\ep})\mathcal{R}:(\Grad_x \vu_{\ep}+\Grad_x^t \vu_{\ep}).
\]
Likewise,
\[
\int_0^T\int_{\Om} \psi \varphi
\mathbb{S}(\vt,\Grad_x\vu):\mathcal{R}(1_{\Om}\vr)
\dxdt
\]
\[
=
\int_0^T\int_{\Om} \psi \varphi
\mu(\vt)(\Grad_x \vu+\Grad_x^t \vu-\Div_x \vu\mathbb{I} ): \mathcal{R}(1_{\Om} \vr)
\dxdt
\]
\[=
\int_0^T\int_{\Om} \psi \mathcal{R}:\left[
\varphi\mu(\vt)(\Grad_x \vu+\Grad_x^t \vu)
\right] \vr
\dxdt
-
\int_0^T\int_{\Om} \psi \varphi
\mu(\vt) \vr\Div_x \vu\dxdt
\]
\beq\label{exi-129}
=
\int_0^T\int_{\Om} \psi \varphi
\mu(\vt) \vr\Div_x \vu\dxdt
+\int_0^T\int_{\Om} \psi \vr \omega(\vt,\vu)\dxdt.
\eeq
Taking advantage of a version of commutator lemma, for instance see Theorem 10.28 \cite{FN1}, it follows that
\beq\label{exi-130}
\|\omega(\vt_{\ep},\vu_{\ep})\|_{ L^1(0,T;W^{p,q}(\Om))  } \leq C \text{   for some   }0<p<1,q>1,
\eeq
uniformly in $\ep$. Application of Div-Curl lemma again, recalled in Section \ref{dicu}, in the time-space domain $(0,T)\times \Om$ with
\[
\mathbf{f}_{\ep}:= (\vr_{\ep} ,\vr_{\ep}\vu_{\ep}),\,\,
\mathbf{g}_{\ep}:= (\omega(\vt_{\ep},\vu_{\ep}),0,0)
\]
gives rise to
\[
\omega(\vt_{\ep},\vu_{\ep}) \vr_{\ep} \rightarrow \overline{\omega(\vt,\vu)}\,\vr \text{   weakly in   }L^1((0,T)\times \Om).
\]
On the other hand, we know from (\ref{exi-80}) and (\ref{exi-105}) that
\beq\label{exi-131}
\overline{\omega(\vt,\vu)}=\omega(\vt,\vu).
\eeq
The identity (\ref{exi-121}) then follows immediately from (\ref{exi-127}), (\ref{exi-128}), (\ref{exi-129}) and (\ref{exi-131}). This finishes the proof of Lemma \ref{effe}.  $\Box$

Recalling (\ref{exi-105}), we may rewrite (\ref{exi-121}) in a more convenient form
\beq\label{exi-132}
\overline{\overline{\mathcal{P}_{M,\delta}(\vr,\vt,\zeta)\vr}}
-\overline{\overline{\mathcal{P}_{M,\delta}(\vr,\vt,\zeta)}}\,\vr
=\mu(\vt)\Big(
\overline{\vr\Div_x \vu}-\vr \Div_x \vu
\Big)
\eeq
a.e. in $(0,T)\times \Om$. It should be noticed that $\mathcal{P}_{M,\delta}(\vr,\vt,\zeta)$ is monotonic increasing with respect to $\vr$. Thus,
\[
\overline{\overline{\mathcal{P}_{M,\delta}(\vr,\vt,\zeta)\vr}}
\geq
\overline{\overline{\mathcal{P}_{M,\delta}(\vr,\vt,\zeta)}}\,\vr;
\]
whence
\beq\label{exi-133}
\overline{\vr\Div_x \vu}\geq \vr \Div_x \vu \text{    a.e. in   }(0,T)\times \Om.
\eeq
Multiplying (\ref{exi-49})$_1$ by $1+\log \vr_{\ep}$, integrating by parts and passing to the limit $\ep\rightarrow 0$ shows
\beq\label{exi-134}
\int_{\Om} \overline{\vr\log \vr}(\tau)\dx
+\int_0^{\tau}\int_{\Om}\overline{\vr\Div_x \vu} \dxdt\leq
\int_{\Om} \vr_{0,\delta}\log \vr_{0,\delta} \dx.
\eeq
On the other hand, it follows from the renormalization theory of transport equations (see \cite{DL}) that the limit continuity equation (\ref{exi-91}) is satisfied in the sense of renormalized solutions. As a consequence,
\beq\label{exi-135}
\int_{\Om} \vr\log \vr(\tau)\dx
+\int_0^{\tau}\int_{\Om}\vr\Div_x \vu \dxdt=
\int_{\Om} \vr_{0,\delta}\log \vr_{0,\delta} \dx .
\eeq
We then conclude from (\ref{exi-134}), (\ref{exi-135}) and the strict convexity of the mapping $\vr \mapsto \vr \log \vr, \,\vr>0$ that
\[
\overline{\vr\log \vr}=\vr\log \vr,
\]
yielding, making use of its convexity again,
\beq\label{exi-136}
\vr_{\ep}\rightarrow \vr \text{    a.e. in    } (0,T)\times \Om.
\eeq

With the help of (\ref{exi-136}), the approximate momentum equation (\ref{exi-120}) is finally reformulated as
\[
\int_0^{T}\int_{\Om}
\left[  \vr\vu \cdot \p_t \vc{\phi}+\vr \vu\otimes \vu :\Grad_x \vc{\phi}+
\left( p_{M,\delta}(\vr,\vt,b)+\f{a}{3}\vt^4 \right)  \Div_x \vc{\phi} \right] \dxdt
\]
\beq\label{exi-137}
=\int_0^{T}\int_{\Om}
\mathbb{S}(\vt,\Grad_x \vu):\Grad_x \vc{\phi} \dxdt
-\int_{\Om} \vr_{0,\delta}\vu_{0,\delta}\cdot \vc{\phi} (0,\cdot)    \dx,
\eeq
for any $\vc{\phi} \in C_c^{\infty}([0,T)\times \Om;\R^2)$.

In addition, by invoking (\ref{exi-117}), (\ref{exi-105}), (\ref{exi-136}), we let $\ep\rightarrow 0$ in (\ref{exi-71}) to obtain the conservation of approximate total energy$\footnote{With (\ref{exi-117}) and (\ref{exi-136}) at hand, pointwise convergence of $\{b_{\ep}\}_{\ep>0}$ is verified immediately. Indeed, it suffices to notice that
\[
|b_{\ep}-b|=|\vr_{\ep}\zeta_{\ep}-\vr \zeta|=|\vr_{\ep}\zeta_{\ep}-\vr_{\ep}\zeta+\vr_{\ep}\zeta-\vr \zeta|\leq \vr_{\ep}|\zeta_{\ep}-\zeta|+\zeta|\vr_{\ep}-\vr|.
\] }$:
\[
\int_{\Om} \left[\f{1}{2}\vr|\vu|^2+\vr e(\vr,\vt) +\f{1}{2}b^2+
\delta \left( \f{\vr^{\Gamma}  }{\Gamma-1}+\vr^2+ \f{b^{\Gamma}  }{\Gamma-1}+b^2 \right)\right](\tau,x)  \dx
\]
\[
=
\int_{\Om} \left[   \f{1}{2}\vr_{0,\delta}|\vu_{0,\delta}|^2+\vr_{0,\delta} e(\vr_{0,\delta},\vt_{0,\delta}) +\f{1}{2}b_{0,\delta}^2
+\delta \left( \f{\vr_{0,\delta}^{\Gamma}  }{\Gamma-1}+\vr_{0,\delta}^2+ \f{b_{0,\delta}^{\Gamma}  }{\Gamma-1}+b_{0,\delta}^2 \right)
\right] \dx
\]
\begin{equation}\label{exi-138}
+\int_0^{\tau}\int_{\Om}
\delta \f{1}{\vt^2}\dxdt
\eeq
for a.e. $\tau\in (0,T)$. By the same  token, the balance of approximate entropy (\ref{exi-109}) is equal to
\[
\int_0^T \int_{\Om}  \Big(\vr s(\vr,\vt)\p_t \phi
+\vr s(\vr,\vt)\vu \cdot \Grad_x \phi \Big)\dxdt
\]
\[
-
\int_0^T \int_{\Om}
\left[\f{\kappa(\vt )}{\vt }+\delta \left(\vt ^{\Gamma-1}+\f{1}{\vt ^2}\right)\right] \Grad_x \vt
\cdot \Grad_x \phi \dxdt
\]
\begin{equation}\label{exi-139}
+\left<\Sigma^{(\delta)} ;\phi\right> = -\int_{\Om} \vr_{0,\delta}s(\vr_{0,\delta},\vt_{0,\delta})\phi(0,\cdot) \dx
\eeq
for any $\phi \in C_c^{\infty}([0,T)\times \overline{\Om}),\phi\geq 0$ and $\Sigma^{(\delta)}$ obeys (\ref{exi-110}).

\subsection{Vanishing artificial pressure}\label{vanpr}
At this level of approximation, we pass to the limit $\delta\rightarrow 0$, obtaining global weak solutions to the original problem (\ref{fmhd_8})-(\ref{fmhd_14}). Since this process is quite similar to the previous one, we shall only present the sketch and point out the main differences.
\subsubsection{Uniform-in-$\delta$ estimates}\label{undel}
From Section \ref{vani} we know that there exists $\{(\vr_{\delta},\vu_{\delta},b_{\delta},\vt_{\delta})\}_{\delta>0}$ solving the approximate equations of density (\ref{exi-91}), magnetic field (\ref{exi-92}), entropy (\ref{exi-139}), momentum (\ref{exi-137}) and total energy (\ref{exi-138}). In this subsection, we denote by $C $ generic positive constants independent of $\delta$. As a consequence of (\ref{exi-8}), (\ref{exi-72})-(\ref{exi-74}), (\ref{exi-82}),
\begin{equation}\label{exi-140}
\sup_{\tau\in (0,T)}
\int_{\Om} \left[\f{1}{2}\vr_{\delta}|\vu_{\delta}|^2+\mathcal{H}_{\overline{\vt}}(\vr_{\delta},\vt_{\delta})+
\f{1}{2}b_{\delta}^2
+
\delta \left( \f{\vr_{\delta}^{\Gamma}  }{\Gamma-1}+\vr_{\delta}^2+ \f{b_{\delta}^{\Gamma}  }{\Gamma-1}+b_{\delta}^2 \right)\right](\tau,x)  \dx\leq C,
\eeq
\begin{equation}\label{exi-141}
\int_0^{T}\int_{\Om}
\f{1}{\vt_{\delta}}\left\{
\mathbb{S}(\vt_{\delta},\Grad_x \vu_{\delta}):\Grad_x \vu_{\delta}
+ \left[\f{\kappa(\vt_{\delta})}{\vt_{\delta}}+\f{\delta }{2} \left(\vt_{\delta}^{\Gamma-1}+\f{1}{\vt_{\delta}^2}\right)\right]|\Grad_x \vt_{\delta}|^2
\right\} \dxdt \leq C,
\eeq
\begin{equation}\label{exi-142}
\int_0^{T}\int_{\Om}
\delta \f{1}{\vt_{\delta}^3}
\dxdt\leq C,
\eeq
\begin{equation}\label{exi-143}
C_{\ast}\vr_{\delta}(t,x) \leq b_{\delta}(t,x) \leq C^{\ast}\vr_{\delta}(t,x) \text{    for a.e.    }(t,x)\in (0,T)\times \Om.
\eeq
In particular, there holds uniform-in-$\delta$ estimates:
\begin{equation}\label{exi-144}
\left\{\begin{aligned}
& \sup_{\tau\in(0,T)}\|\sqrt{\vr_{\delta}}\,\vu_{\delta}(\tau)\|_{L^2(\Om)}+
\sup_{\tau\in(0,T)}\|\vt_{\delta}(\tau)\|_{L^{4}(\Om)}\leq C, \\
&\sup_{\tau\in(0,T)}\|\vr_{\delta}(\tau)\|_{L^{\gamma}(\Om)}+
\sup_{\tau\in(0,T)}\|b_{\delta}(\tau)\|_{L^{2}(\Om)}\leq C,\\
&  \sup_{\tau\in(0,T)}\|(\vr_{\delta},b_{\delta})(\tau)\|_{L^{\Gamma}(\Om)}\leq C \delta^{-1/\Gamma  }, \\
& \|\Grad_x \log \vt_{\delta}\|_{L^2((0,T)\times \Om  ) } +\|\Grad_x \vt_{\delta}^{3/2}\|_{L^2((0,T)\times \Om  ) }\leq C,\\
& \|\Grad_x\vu_{\delta}+\Grad_x^t\vu_{\delta}-\Div_x \vu_{\delta}\mathbb{I}  \|_{L^2((0,T)\times \Om  ) }\leq C ;\\
\end{aligned}\right.
\end{equation}
whence, with the help of generalized Korn's inequality (\ref{ap-1}) and generalized Poincar\'{e} inequality (\ref{ap-4}),
\begin{equation}\label{exi-145}
\|\vu_{\delta}\|_{L^2(0,T;W_0^{1,2}(\Om))  }  \leq C     ,
\eeq
\begin{equation}\label{exi-146}
\left\|\left(\vt_{\delta},\vt_{\delta}^{3/2}\right)\right\|_{L^2(0,T;W^{1,2}(\Om))  }  \leq C     ,
\eeq
\begin{equation}\label{exi-147}
\|\mathbb{S}(\vt_{\delta},\Grad_x \vu_{\delta})\|_{  L^2(0,T;  L^{4/3} (\Om)) }\leq C,
\eeq
\begin{equation}\label{exi-148}
\|\log \vt_{\delta}\|_{L^2(0,T;W^{1,2}(\Om))  }  \leq C .
\eeq
Furthermore, in analogy with (\ref{exi-111}), we may choose
\[
\vc{\phi}(t,x)=\psi(t)\mathcal{B}\left( \vr_{\delta}^{\beta}-\f{1}{|\Om|} \int_{\Om}\vr_{\delta}^{\beta} \dx  \right),\,\,
\psi(t)\in C_c^{\infty}( (0,T) )
\]
as a test function in (\ref{exi-137}) to find that
\begin{equation}\label{exi-149}
\int_0^T \psi \int_{\Om}\left( p(\vr_{\delta},\vt_{\delta})+\f{1}{2}b_{\delta}^2 +\delta \left(\vr_{\delta}^{\Gamma}+\vr_{\delta}^2+b_{\delta}^{\Gamma}+b_{\delta}^2\right)
\right) \vr_{\delta}^{\beta} \dxdt \leq C
\eeq
for some $\beta>0$. The details are again omitted here for simplicity, see \cite{FNP,FN1}.

Similar to (\ref{exi-78})-(\ref{exi-90}), we infer from (\ref{exi-140})-(\ref{exi-149}) that there exists $(\vr,\vu,b,\vt)$ such that as $\delta  \rightarrow 0$,
\begin{equation}\label{exi-150}
\vr_{\delta} \rightarrow \vr \text{    weakly    }-\ast \text{  in }L^{\infty}(0,T;L^{\gamma}(\Om)), \,\,
\text{   weakly in   }L^{ \max\{\gamma+\beta,2+\beta  \}} ((0,T)\times \Om),
\eeq
\begin{equation}\label{exi-151}
b_{\delta} \rightarrow b \text{     weakly   }-\ast \text{  in }L^{\infty}(0,T;L^{2}(\Om)),
\,\,
\text{   weakly in   }L^{ \max\{\gamma+\beta,2+\beta  \}} ((0,T)\times \Om),
\eeq
\begin{equation}\label{exi-152}
\vu_{\delta} \rightarrow \vu \text{   weakly in   }L^2(0,T;W_0^{1,2}(\Om;\R^2)),
\eeq
\begin{equation}\label{exi-153}
\vt_{\delta} \rightarrow \vt \text{  weakly}-\ast \text{  in }L^{\infty}(0,T;L^{4}(\Om)),
\eeq
\begin{equation}\label{exi-154}
C_{\ast}\vr(t,x) \leq b(t,x) \leq C^{\ast}\vr(t,x) \text{    for a.e.    }(t,x)\in (0,T)\times \Om,
\eeq
\begin{equation}\label{exi-155}
\vr_{\delta} \vu_{\delta} \rightarrow \vr \vu \text{    weakly    }-\ast \text{  in }L^{\infty}(0,T;L^{ \f{2\max\{ 2,\gamma\}}{ \max\{ 2,\gamma\} +1}  }(\Om)),
\eeq
\begin{equation}\label{exi-156}
\vr_{\delta} \vu_{\delta}\otimes \vu_{\delta}  \rightarrow \vr \vu\otimes \vu \text{    weakly in   } L^p((0,T)\times \Om), \text{   for some } p>1,
\eeq
\begin{equation}\label{exi-157}
b_{\delta} \vu_{\delta} \rightarrow b \vu \text{    weakly in   } L^p((0,T)\times \Om), \text{   for some } p>1,
\eeq
\begin{equation}\label{exi-158}
\|\delta(\vr_{\delta}^{\Gamma},\vr_{\delta}^{2},b_{\delta}^{\Gamma},b_{\delta}^{2},\vt_{\delta}^{-2})  \|_{L^1((0,T)\times \Om)} \rightarrow 0 .
\eeq

Therefore, we pass to the limits $\delta \rightarrow 0 $ in (\ref{exi-91})-(\ref{exi-92}), obtaining the equations of continuity and magnetic field:
\beq\label{exi-159}
\int_0^{T}\int_{\Om} \Big( \vr \p_t \phi +\vr \vu \cdot \Grad_x \phi\Big)\dxdt
+  \int_{\Om} \vr_{0,\delta}\phi(0,\cdot) \dx = 0,
\eeq
\beq\label{exi-160}
\int_0^{T}\int_{\Om} \Big( b \p_t \phi +b \vu \cdot \Grad_x \phi \Big)\dxdt
+  \int_{\Om} b_{0,\delta}\phi(0,\cdot) \dx = 0,
\eeq
for any $\phi \in C_c^{\infty}( [0,T)\times \overline{\Om}  )$. Next, letting $\delta \rightarrow 0 $ in (\ref{exi-137}) gives the balance of momentum:
\[
\int_0^{T}\int_{\Om}
\left[  \vr\vu \cdot \p_t \vc{\phi}+\vr \vu\otimes \vu :\Grad_x \vc{\phi}+
\left( \overline{p_{M}(\vr,\vt,b)}+\f{a}{3}\overline{\vt^4} \right)  \Div_x \vc{\phi} \right] \dxdt
\]
\beq\label{exi-161}
=\int_0^{T}\int_{\Om}
\overline{\mathbb{S}(\vt,\Grad_x \vu)}:\Grad_x \vc{\phi} \dxdt
-\int_{\Om} \vr_{0,\delta}\vu_{0,\delta}\cdot \vc{\phi} (0,\cdot)    \dx,
\eeq
for any $\vc{\phi} \in C_c^{\infty}([0,T)\times \Om;\R^2)$. To proceed, we pass to the limit $\delta \rightarrow 0 $ in (\ref{exi-138}), attaining the balance of total energy:
\[
\int_{\Om} \left(\f{1}{2}\vr|\vu|^2+\overline{\vr e(\vr,\vt)} +\f{1}{2}\overline{b^2} \right)
(\tau,x)  \dx
\]
\begin{equation}\label{exi-162}
=
\int_{\Om} \left(   \f{1}{2}\vr_{0}|\vu_{0}|^2+\vr_{0} e(\vr_{0},\vt_{0}) +\f{1}{2}b_{0}^2
\right) \dx,
\eeq
for a.e. $\tau\in (0,T)$.

\subsubsection{Strong convergence of $\{\vt_{\delta}\}_{\delta>0}$}\label{strovt-2}
In order to pass to the limit in the approximate entropy equation (\ref{exi-139}), we notice first that (\ref{exi-146}) and Sobolev's embedding inequality show
\begin{equation}\label{exi-163}
\|\vt_{\delta}\|_{L^3(0,T; L^{\f{3}{2}q} (\Om))   }\leq C, \text{  for any }1\leq q<\infty.
\eeq
Combining (\ref{exi-163}), the structural hypothesis (\ref{fmhd_21}) with the uniform-in-$\delta$ estimates obtained so far,
\begin{equation}\label{exi-164}
\| (\vr_{\delta}s(\vr_{\delta},\vt_{\delta}),\vr_{\delta}s(\vr_{\delta},\vt_{\delta})\vu_{\delta}  )  \|_
{L^p((0,T)\times \Om )} \leq C, \text{   for some   } p>1,
\eeq
\begin{equation}\label{exi-165}
\left\| \f{\kappa(\vt_{\delta})}{\vt_{\delta}}\Grad_x \vt_{\delta} \right \|_
{L^p((0,T)\times \Om )} \leq C, \text{   for some   } p>1.
\eeq
Moreover, in view of (\ref{exi-141}),
\begin{equation}\label{exi-166}
\|\delta(\vt_{\delta}^{\Gamma-1},\vt_{\delta}^{-2}    ) \Grad_x \vt_{\delta}  \|_{L^1((0,T)\times \Om)} \rightarrow 0.
\eeq
With these estimates, one then follows the same line as in Section \ref{strovt}, by Div-Curl lemma and parameterized Young measures, to conclude
\begin{equation}\label{exi-167}
\vt_{\delta} \rightarrow \vt \text{     a.e. in   } (0,T)\times \Om.
\eeq

Taking into account of (\ref{exi-150})-(\ref{exi-153}), (\ref{exi-164})-(\ref{exi-167}) and dropping the $\delta$-dependent terms in $\Sigma^{(\delta)}$, we let $\delta \rightarrow 0 $ in (\ref{exi-139}) to see that
\[
\int_0^{T}\int_{\Om}
\left(
\overline{\vr s(\vr,\vt) } \,\p_t \phi +\overline{\vr s(\vr,\vt)}\, \vu \cdot \Grad_x \phi +\f{\mathbf{q}\cdot \Grad_x \phi}{\vt}
\right) \dxdt
\]
\beq\label{exi-168}
+\int_0^{T}\int_{\Om}
\f{\phi}{\vt}\left(
\mathbb{S}(\vt,\Grad_x \vu):\Grad_x \vu- \f{\mathbf{q} \cdot \Grad_x \vt}{\vt}
\right)\dxdt
\leq
\int_{\Om}
\vr_0 s(\vr_0,\vt_0) \phi  (0,\cdot)\dx
\eeq
for any $\phi \in C^{\infty}([0,T)\times \overline{\Om})$, $\phi \geq 0$.

\subsubsection{Strong convergence of $\{\vr_{\delta}\}_{\delta>0}$}\label{strovr-2}
Clearly, it remains to show pointwise convergences of $\{(\vr_{\delta},b_{\delta})\}_{\delta>0}$. Basically, this is finished as before with the aid of variable reduction$\footnote{Compared with the compressible Navier-Stokes system \cite{FNP}, full compressible two-fluid model \cite{KNC}, our MHD model naturally gives better integrability of the approximate densities $\{\vr_{\delta}\}_{\delta>0}$, due to the nice integrability property of appoximate magnetic fields $\{b_{\delta}\}_{\delta>0}$ and the domination condition (\ref{exi-143}). }$. More specifically, we first set

\[
p_{M}(\vr,\vt,b):=p_M(\vr,\vt)+\f{1}{2}b^2
\]
\[
=\vr^{\gamma}+\vr\vt+\f{1}{2}b^2 ,
\]

\[
\zeta_{\delta}:=
\begin{cases}
\f{b_{\delta}}{\vr_{\delta}}& \text{ if } \vr_{\delta}>0, \\
\f{C_{\ast}+C^{\ast}   }{2}& \text{ if }\vr_{\delta}=0,\\
\end{cases}
\,\,\,\,
\zeta:=
\begin{cases}
\f{b}{\vr}& \text{ if } \vr>0, \\
\f{C_{\ast}+C^{\ast}   }{2}& \text{ if }\vr=0,\\
\end{cases}
\]
and rewrite $p_{M}(\vr,\vt,b)$ in an equivalent form
\[
\mathcal{P}_{M}(\vr,\vt,\zeta)=\vr^{\gamma}+\vr\vt+\f{1}{2}\zeta^2\vr^2.
\]
As a consequence,
\[
p_{M}(\vr_{\delta},\vt_{\delta},b_{\delta})
=\mathcal{P}_{M}(\vr_{\delta},\vt_{\delta},\zeta_{\delta})
\]
\beq\label{exi-169}
=\mathcal{P}_{M}(\vr_{\delta},\vt_{\delta},\zeta_{\delta})
-\mathcal{P}_{M}(\vr_{\delta},\vt,\zeta)+\mathcal{P}_{M}(\vr_{\delta},\vt,\zeta).
\eeq
Similar to Lemma \ref{comp-lem}, it still holds the following crucial estimate at this level
\beq\label{exi-170}
\int_{\Om}\vr_{\delta} |\zeta_{\delta}-\zeta|^p (\tau,\cdot)\dx \rightarrow 0
\eeq
for any $\tau\in [0,T]$ and any $1\leq p<\infty$; whence
\[
\int_0^T\int_{\Om}\left|
\mathcal{P}_{M}(\vr_{\delta},\vt_{\delta},\zeta_{\delta})
-\mathcal{P}_{M}(\vr_{\delta},\vt,\zeta)\right|
\dxdt
\]
\[
\leq C \left(
\int_0^T\int_{\Om} \vr_{\delta}|\vt_{\delta}-\vt| \dxdt
+\int_0^T\int_{\Om}
\vr_{\delta}^2
 |\zeta_{\delta}-\zeta| \dxdt
\right)
\]
\beq\label{exi-171}
 \rightarrow 0 \,\,\text{     as         }\delta \rightarrow 0,
\eeq
due to (\ref{exi-150}) and (\ref{exi-167}). Thus, the balance of momentum (\ref{exi-161}) may be reformulated as
\[
\int_0^{T}\int_{\Om}
\left[  \vr\vu \cdot \p_t \vc{\phi}+\vr \vu\otimes \vu :\Grad_x \vc{\phi}+
\left( \overline{\overline{\mathcal{P}_{M}(\vr,\vt,\zeta)}}+\f{a}{3}\vt^4 \right)  \Div_x \vc{\phi} \right] \dxdt
\]
\beq\label{exi-172}
=\int_0^{T}\int_{\Om}
\mathbb{S}(\vt,\Grad_x \vu):\Grad_x \vc{\phi} \dxdt
-\int_{\Om} \vr_{0,\delta}\vu_{0,\delta}\cdot \vc{\phi} (0,\cdot)    \dx,
\eeq
for any $\vc{\phi} \in C_c^{\infty}([0,T)\times \Om;\R^2)$.

To proceed, let $\{T_k(z)\}_{k\geq 1}$ be a sequence of cut-off functions
\[
T_k(z)=kT\left( \f{z}{k} \right),
\]
for some smooth concave function $T(z)$ on $[0,\infty)$ obeying
\[
T(z)=
\begin{cases}
z& \text{ if } 0\leq z\leq 1, \\
2& \text{ if }z\geq 3.\\
\end{cases}
\]
Following step by step the proof of Lemma \ref{effe}, with slight modifications, we also obtain the ``effective viscous flux identity" at this level:
\beq\label{exi-173}
\overline{\overline{\mathcal{P}_{M}(\vr,\vt,\zeta)T_k(\vr)}}
-\overline{\overline{\mathcal{P}_{M}(\vr,\vt,\zeta)}}\,\overline{T_k(\vr)}
=\mu(\vt)\Big(
\overline{T_k(\vr)\Div_x \vu}-\overline{T_k(\vr)}\, \Div_x \vu
\Big)
\eeq
a.e. in $(0,T)\times \Om$. In view of (\ref{exi-145}) and (\ref{exi-150}),
\beq\label{exi-174}
\vr_{\delta} \in L^{ \max\{\gamma+\beta,2+\beta  \}} ((0,T)\times \Om),\,\, \vu_{\delta}\in L^2(0,T;W_0^{1,2}(\Om;\R^2)).
\eeq
Based on (\ref{exi-174}), we apply the renormalization theory of transport equation due to DiPerna and Lions (see \cite{DL}) to the approximate continuity equation (\ref{exi-91}) to infer that
\[
\int_0^T\int_{\Om} \Big(\overline{  \vr L_k(\vr)} \, \p_t \phi+ \overline{  \vr L_k(\vr)} \, \vu\cdot \Grad_x \phi - \overline{  T_k(\vr)\Div_x \vu} \,\phi  \Big)\dxdt
\]
\beq\label{exi-175}
=-\int_{\Om}  \vr_0 L_k(\vr_0)\phi(0,\cdot)\dx;
\eeq
similarly,
\[
\int_0^T\int_{\Om} \Big( \vr L_k(\vr) \, \p_t \phi+  \vr L_k(\vr)\, \vu\cdot \Grad_x \phi -   T_k(\vr)\Div_x \vu \,\phi  \Big)\dxdt
\]
\beq\label{exi-176}
=-\int_{\Om}  \vr_0 L_k(\vr_0)\phi(0,\cdot)\dx,
\eeq
for any $\phi \in C_c^{\infty}( [0,T)\times \overline{\Om} )$ and $L_k(\vr)$ is defined by
\[
L_k(\vr)=\int_1^{\vr}\f{T_k(z)}{z^2}\text{d}z.
\]
Combining (\ref{exi-175}) and (\ref{exi-176}), we arrive at
\[
\int_{\Om} \left(\overline{  \vr L_k(\vr)} -\vr L_k(\vr) \right)(\tau)\dx
+\int_0^{\tau} \int_{\Om} \left( \overline{  T_k(\vr)\Div_x \vu}- \overline{T_k(\vr)}\,\Div_x \vu \right)
\dxdt
\]
\beq\label{exi-177}
=\int_0^{\tau} \int_{\Om} \left(
T_k(\vr)\Div_x \vu- \overline{T_k(\vr)}\,\Div_x \vu \right)\dxdt
\eeq
for a.e. $\tau\in (0,T)$. On the one hand, by (\ref{exi-167}), (\ref{exi-173}) and the monotonic increasing of $\mathcal{P}_{M}(\vr,\vt,\zeta)$ with respect to $\vr$, it holds
\beq\label{exi-178}
\overline{T_k(\vr)\Div_x \vu}-\overline{T_k(\vr)}\, \Div_x \vu \geq 0
\eeq
a.e. in $(0,T)\times \Om$. On the other hand, it is easy to check, using the definition of $T_k(\cdot)$, that as $k\rightarrow \infty$
\beq\label{exi-179}
T_k(\vr)\rightarrow \vr,\,\, \overline{T_k(\vr)}\rightarrow \vr \text{    strongly in   }L^1( (0,T)\times \Om  ).
\eeq
Observe next that
\[
\|T_k(\vr)- \overline{T_k(\vr)}\| _{ L^{ \max\{\gamma+\beta,2+\beta  \}} ((0,T)\times \Om) }
\]
\[
\leq \lim_{\delta\rightarrow 0}\inf \| T_k(\vr)-T_k(\vr_{\delta})\|_{ L^{ \max\{\gamma+\beta,2+\beta  \}} ((0,T)\times \Om) }
\]
\beq\label{exi-180}
\leq C \text{    uniformly in   }k,\delta,
\eeq
in accordance with (\ref{exi-150}). We then deduce from (\ref{exi-152}), (\ref{exi-179}) and (\ref{exi-180}) that
\[
\int_0^{\tau} \int_{\Om} \left(
T_k(\vr)\Div_x \vu- \overline{T_k(\vr)}\,\Div_x \vu \right)\dxdt
\rightarrow 0 \text{    as    }k\rightarrow \infty,
\]
which together with (\ref{exi-177}), (\ref{exi-178}) yields, upon passing $k\rightarrow \infty$,
\beq\label{exi-181}
\vr_{\delta}\rightarrow\vr \text{     a.e. in  } (0,T)\times \Om.
\eeq
In particular, pointwise convergence of $\{b_{\delta}\}_{\delta>0}$ is verified as in the footnote of equation (\ref{exi-138}). The proof of our main Theorem \ref{ls_1} is finished completely.

\section{Concluding remarks}\label{conre}

\begin{itemize}
\item{
The structural hypotheses imposed in Section \ref{struct} is just to simplify the presentation. Indeed, the strategy adopted in this paper allows us to treat more general case, such as:
\[
\mu=\mu(\vt)\sim (1+\vt^{\alpha}), \,\,0\leq \eta=\eta(\vt)\leq \eta_1 (1+\vt^{\alpha}),
\]
\[
\kappa=\kappa(\vt)=\kappa_M(\vt)+\kappa_R(\vt),\,\, \kappa_M(\vt)\sim (1+\vt^{\alpha}),
\,\,\kappa_R(\vt)\sim (1+\vt^{3}),
\]
with suitable $\alpha \in (0,1]$.
}
\item{
In order to prove the \emph{existence} of global weak solutions with finite energy initial data to multi-dimensional compressible, viscous, non-resistive full MHD system, we consider the case where the motion of fluids takes place in the plane, while the magnetic field acts on the fluids only in the vertical direction. Under this symmetry, the magnetic equation without resistivity reduces to the continuity equation, transported by the same velocity field like the density. Then we are allowed to invoke the techniques from compressible two-fluid model and full Navier-Stokes system so as to construct weak solution through three-level approximations. However, for the \emph{general} 2D/3D compressible, viscous, non-resistive isentropic/full MHD system, the existence of global weak solutions with large initial data remains completely open.
}

\end{itemize}

\centerline{\bf Acknowledgements}
The research of Y. Sun is supported by the NSF of China under grant numbers 11571167, 11771395, 11771206 and PAPD of Jiangsu Higher Education Institutions.

\centerline{\bf Conflicts of interest}
The authors declare that they have no conflicts of interest.

\section{Appendix}\label{appe}
For convenience of the reader, we list some necessary mathematical tools in this appendix without giving the detailed proofs.

\subsection{Helmholtz function}
Let $\overline{\vr},\overline{\vt}>0$. The Helmholtz function is defined as
\[
\mathcal{H}_{\overline{\vt}}(\vr,\vt):=\vr e(\vr,\vt)-\overline{\vt}\vr s(\vr,\vt).
\]
There holds:
\begin{itemize}
\item{
the mapping $\vr \mapsto \mathcal{H}_{\overline{\vt}}(\vr,\overline{\vt})$ is strictly convex for any $\vr>0$;
}
\item{
the mapping $\vt \mapsto \mathcal{H}_{\overline{\vt}}(\vr,\vt)$ is decreasing on $(0,\overline{\vt})$, while increasing on $(\overline{\vt},\infty)$;
}
\item{ coercivity property:
\beq\label{ap-3}
\mathcal{H}_{\overline{\vt}}(\vr,\vt)
\geq \f{1}{4} \left(
\vr e(\vr,\vt) +\overline{\vt} \vr |s(\vr,\vt)|
\right)
-\left|
(\vr-\overline{\vr})\f{\p \mathcal{H}_{  2\overline{\vt}  } }
{\p \vr}(\overline{\vr},2\overline{\vt})
+\mathcal{H}_{2\overline{\vt}}(\overline{\vr},2\overline{\vt})
\right|
\eeq
for any $\vr,\vt>0$.
}
\end{itemize}

The above properties of Helmholtz functions are proved essentially based on Gibbs' relation (\ref{fmhd_6_1}) and thermodynamic stability conditions (\ref{fmhd_6_2}). The interested reader may consult \cite{FN1} for the details.

\subsection{Generalized Poincar\'{e} inequality}
Let $N\geq 2$ be an integer and $1\leq p\leq \infty$ and $\Om\subset \R^N$ be a bounded Lipschitz domain. Suppose also that $1<\tilde{p}<\infty$ and $\Om' \subset \Om$ is measurable with $|\Om'|>0$. Then there exists $C>0$ depending only on $p,\tilde{p},\Om'$ such that
\beq\label{ap-4}
\|u\|_{W^{1,p}(\Om)}
\leq C \left(
\|\Grad_x u\|_{L^p(\Om)}
+\|u\|_{L^{\tilde{p}}(\Om')}
\right)
\eeq
for any $u\in W^{1,p}(\Om)$.

We refer to Theorem 10.14 in Ref. \cite{FN1} for the proof.

\subsection{Generalized Korn's inequality}

Let $N\geq 2$ be an integer and $1<p<\infty$. Then there exists $C>0$ depending only on $p,N$ such that
\beq\label{ap-1}
\|\Grad_x \mathbf{U} \|_{L^p(\R^N)}
\leq C
\left\|\Grad_x \mathbf{U}+\Grad_x^t \mathbf{U}-\f{2}{N}\Div_x \mathbf{U} \mathbb{I}  \right\|_{L^p(\R^N)}
\eeq
for any $\mathbf{U}\in W^{1,p}(\R^N;\R^N)$. Moreover, assuming that $\Om\subset \R^N$ is a bounded Lipschitz domain, it also holds for some $C>0$ depending only on $p,N,\Om$ that
\beq\label{ap-2}
\|\Grad_x \mathbf{U} \|_{L^p(\Om)}
\leq C
\left\|\Grad_x \mathbf{U}+\Grad_x^t \mathbf{U}-\f{2}{N}\Div_x \mathbf{U} \mathbb{I}  \right\|_{L^p(\Om)}
\eeq
for any $\mathbf{U}\in W_0^{1,p}(\Om;\R^N)$.

The proof of (\ref{ap-1}) may be found in Theorem 10.16 in Ref. \cite{FN1}; whence (\ref{ap-2}) follows immediately.

\subsection{Div-Curl lemma}\label{dicu}
Let $N\geq 2$ be an integer and $\Om \subset \R^N$ be open. Suppose that
\[
\mathbf{f}_n \rightarrow \mathbf{f} \text{   weakly in   } L^p(\Om;\R^N),
\]
\[
\mathbf{g}_n \rightarrow \mathbf{g} \text{   weakly in   } L^q(\Om;\R^N),
\]
with
\[
\f{1}{p}+\f{1}{q}=\f{1}{r}<1.
\]
Suppose also that
\[
\{\Div_x \mathbf{f}_n\}_{n \geq 1},\,\, \{\mathbf{curl}_x \mathbf{g}_n\}_{n \geq 1}
\]
are precompact in $W^{-1,\tilde{p}}(\Om)$ for some $\tilde{p}>1$. Then we have
\[
\mathbf{f}_n \cdot \mathbf{g}_n \rightarrow \mathbf{f} \cdot \mathbf{g} \text{   weakly in   } L^r(\Om).
\]

We refer to Theorem 10.21 in Ref. \cite{FN1} for the proof.



\end{document}